\RequirePackage{fix-cm}
\documentclass[smallextended]{svjour3}       
\smartqed  
\usepackage{lineno,hyperref}
\usepackage{graphicx}
\usepackage{csquotes}
\usepackage{epsfig}
\usepackage{lscape}
\usepackage{amsfonts}
\usepackage[misc,geometry]{ifsym}
\usepackage{color}
\usepackage{amsmath,amssymb,dsfont}
\usepackage{listings}
\usepackage{rotating}
\usepackage{multicol}
\usepackage[english]{babel}
\usepackage{lineno}
\usepackage{enumerate}
\usepackage[justification=centering]{caption}
\usepackage{caption}
\usepackage{subcaption}
\modulolinenumbers[5]
\newcommand{\be}{\begin{equation}}
\newcommand{\ee}{\end{equation}}
\newcommand{\bea}{\begin{eqnarray}}
\newcommand{\eea}{\end{eqnarray}}
\newcommand{\nn}{\nonumber}
\newcommand{\bee}{\begin{eqnarray*}}
\newcommand{\eee}{\end{eqnarray*}}

\newcommand{\bt}{\begin{tabbing}}
\newcommand{\et}{\end{tabbing}}
\newcommand{\btb}{\begin{tabular}}
\newcommand{\etb}{\end{tabular}}

\newcommand{\bc}{\begin{center}}
\newcommand{\ec}{\end{center}}

\newtheorem{t1}{Theorem}
\newcommand{\lb}{\label}

\newtheorem{d1}{Definition}
\newtheorem{l1}{Lemma}

\newtheorem{alg1}{Algorithm}
\begin{document}

\title{A Sequential Quadratic Programming Method for Constrained Multi-objective Optimization Problems}

\author{Md Abu Talhamainuddin Ansary       \and
        Geetanjali Panda 
}


\institute{Md Abu Talhamainuddin Ansary(\Letter) \at
              Department of Mathematics, Indian Institute of Technology Kharagpur, Khragpur, India\\
              \email{md.abutalha2009@gmail.com}
           \and
           Geetanjali Panda   \at
             Department of Mathematics, Indian Institute of Technology Kharagpur, Kharagpur, India\\
\email{geetanjali@maths.iitkgp.ernet.in}
}

\date{Received: date / Accepted: date}

\maketitle

\begin{abstract}
In this article, a globally convergent sequential quadratic programming (SQP) method is developed for multi-objective optimization problems
with inequality type constraints. A feasible descent direction is obtained using a linear approximation of all objective functions as well as constraint functions. The sub-problem at every iteration of the sequence has feasible solution. A non-differentiable penalty function is used to deal with constraint violations. A descent sequence is generated which converges to a critical point under the Mangasarian-Fromovitz constraint qualification along with some other mild assumptions. The method is compared with a selection of existing methods on a suitable set of test problems.
\keywords{multi-objective optimization\and SQP method \and critical point\and Mangasarian-Fromovitz constraint qualification \and purity metric \and spread metrics}
 \subclass{90C26 \and  	49M05 \and 97N40 \and 	90B99}
\end{abstract}

\section{Introduction}
\label{intro}
A widely used line search technique for solving constrained single objective optimization problems is SQP method, which was developed by Wilson in 1963 and modified by several researchers (see \cite{mang1,rob1}) in various directions. A serious limitation of these methods is the inconsistency of the quadratic sub-problem. Powell (\cite{powell}) suggested a modified sub-problem to overcome this restriction, which was further modified in \cite{burke1,liu2,mo1} for better efficiency. SQP method in \cite{burke1} converges to an infeasible point in some situations. SQP method in \cite{liu1} is a two step method. But, SQP method in \cite{mo1} is one step method and always converges to a feasible point. These developments are related to single objective optimization problems. In this article, a convergent SQP iterative scheme is developed for constrained multi-objective optimization problems, in the light of \cite{mo1}.\\

Classical methods for solving multi-objective optimization problems are either scalarization methods or heuristic methods. Scalarization methods reduce the multi-objective optimization problem to a single objective optimization problem using pre determined parameters. Heuristic methods do not guarantee the convergence to the solution. To address these limitations, line search methods for unconstrained multi-objective optimization problems have been developed since 2000 by many researchers (\cite{mat1,mat2,flg1,flg2,pav1}), which are treated as the extension of single objective line search techniques. Possible extension of these concepts to constrained multi-objective problems is an interesting area of research in recent times.\\

The steepest descent method for multi-objective problems developed by Fliege and Svaiter (\cite{flg2}) uses the linear approximation of all objective functions to find a descent direction. This concept is extended in \cite{drummond1} to projected gradient method for vector optimization problems, which is further extended in \cite{fukuda1,cruz1} in different directions. An interior point algorithm is developed in \cite{miglierina1} for box constrained multi-objective optimization problems using the concept of vector pseudo gradient. Recently Fliege and Vaz (\cite{flg3}) and Gebken et al. (\cite{Bennet1}) have developed SQP methods for constrained multi-objective optimization problems using the ideas of single objective SQP methods. The sub-problem in \cite{flg3} is not necessarily feasible at every iteration step. Some restoration process is used to make the sub-problem feasible, and approximate Pareto front is generated. The SQP method in \cite{Bennet1} requires feasible initial approximation, which is very difficult in nonlinear constrained problems. In addition to this, the iterative process in \cite{Bennet1} does not use penalty function, and only active constraints are used in the sub-problem. In this article these difficulties are taken care. A modified SQP scheme is developed using a different sub-problem so that the infeasibility of the sub-problem at every iteration step can be avoided and a non-differentiable penalty function is used to restrict constraint violations.\\

The outline of this article is as follows. Some preliminaries on the existence of the solution of a multi-objective optimization problem are discussed in Section \ref{secpre}. A modified SQP scheme for inequality constrained multi-objective optimization problems is developed in Section \ref{secsqp} and global convergence of this scheme is proved in Section \ref{secgconv}. In Section \ref{secnex}, the proposed method is compared with existing methods using a set of test problems.
\section{Preliminaries}
\lb{secpre}
Consider the multi-objective optimization problem:
\bee
MOP:~~& &~\min~~\left(f_1(x),f_2(x),...,f_m(x)\right) \\
   & & ~~\mbox{s.t. }g_i(x)\leq 0,~~~i=1,2,...,p,
\eee
where $f_j,g_i:\mathbb{R}^n \rightarrow \mathbb{R}$ are continuously differentiable for $j\in\{1,2,...,m\}$ and $i\in\{1,2,...,p\}$. Denote $\mathbb{R}^n_{+}=\{x\in\mathbb{R}^n|x_i\geq 0, i=1,2,...,n\}$,\\$\Lambda_n=\{1,2,...,n\}$ for any $n\in \mathbb{N}$,
and the feasible set of the $MOP$ as\\$X=\{ x\in\mathbb{R}^n:g_i(x)\leq 0,~i\in \Lambda_p\}$. Inequality in $\mathbb{R}^n$ is understood componentwise. If there exists $x\in X$ which minimizes all objective functions simultaneously then it is an ideal solution. But in practice, decrease of one objective function may cause increase of another objective function. So, in the theory of multi-objective optimization, optimality is replaced by efficiency. A feasible point $x^*\in X$ is said to be an efficient solution of the $MOP$ if there does not exist $x\in X$ such that $f(x)\leq f(x^*)$ and $f(x)\neq f(x^*)$
hold where $f(x)=(f_1(x),f_2(x),...,f_m(x))$. A feasible point $x^*\in X$ is said to be a weak efficient solution of $MOP$ if there does not exist $x\in X$ such that $f(x)< f(x^*)$ holds. For $x, y \in X$, we say $y$ dominates $x$, if and only if $f(y)\leq f(x)$, $ f(y)\neq f(x)$. A point $x\in X$ is said to be non dominated if there does not exist any $y\in X$ such that $y$ dominates $x$.  If $X^*$ is the set of all efficient solutions of the $MOP$, then the image of $X^*$ under $f$, i.e. $f(X^*)$ is said to be the Pareto front of the $MOP.$\\

In our analysis, we use the $L_{\infty}$ non-differentiable penalty function
 \bee
  \Phi(x)=\max\{ 0,g_i(x),i\in\Lambda_p\}.\lb{phi}
 \eee
In order to obtain a feasible descent direction, the penalty function for the $MOP$ is used as the following merit function $\Psi_{j,\sigma}(x) $, with a penalty parameter $\sigma>0$ as
\bee
\Psi_{j,\sigma}(x)=f_j(x)+\sigma \Phi(x),~~j\in\Lambda_m.
\eee
Let $I(x)=\{i \in \Lambda_p: g_i(x)=\Phi(x)\}$ be the set of active or most violated constraints.  The directional derivative of $\Phi(x)$  in  any direction $d\in\mathbb{R}^n$ is
\bee
\Phi'(x;d)=\underset{i\in I(x)}{\max}\{\nabla g_i(x)^Td\},
\eee
In general $\Phi'(x;d)$ is not continuous. A continuous approximation of $\Phi'(x;d)$ (see \cite{baz1}) is
\bee
\Phi^*(x;d)=\underset{i\in I(x)}{\max}\{g_i(x)+\nabla g_i(x)^Td,0\}-\Phi(x).
\eee
 Thus, an approximation of the directional derivative of $\Psi_{j,\sigma}(x)$ is
 \bee
\theta_{j,\sigma}(x;d)=\nabla f_j(x)^Td+\sigma \Phi^*(x;d),~~j\in\Lambda_m.\lb{theta}
\eee
If all $f_j,g_i$ are continuously differentiable then the necessary condition of weak efficiency can be derived using  Motzkin's theorem as follows.
\begin{t1}(\textbf{Fritz John necessary condition [Theorem 3.1.1,\cite{kmm1}]})\\
Suppose $f_j$, $j\in\Lambda_m$ and $g_i$, $i\in\Lambda_p$ are continuously differentiable at $x^* \in X$. If $x^*$ is a weak efficient solution of the $MOP$  then there exists $(\lambda,\mu)\in \mathbb{R}_{+}^m\times \mathbb{R}_{+}^p$, $(\lambda,\mu)\neq 0^{m+p}$  satisfying
\bea
 \underset{j \in\Lambda_m}{\Sigma} \lambda_j \nabla f_j(x^*)+\underset{i\in\Lambda_p}{\Sigma} \mu_i \nabla g_i(x^*)&=&0\lb{fj1}\\
 \mu_i g_i(x^*)&=&0 ~~\forall ~~i\in\Lambda_p\lb{fj2}.
\eea
\end{t1}
The set of the vector $(\lambda,\mu)\in \mathbb{R}^{m}_{+} \times \mathbb{R}^{p}_{+}\setminus\{0^{m+p}\}$ satisfying (\ref{fj1}) and (\ref{fj2}) are called Fritz John multipliers associated with $x^*$. But the Fritz John necessary condition does not guarantee $\lambda_j>0$, for at least one $j$. So  some constraint qualifications or regularity conditions should hold to ensure it.\\
Several constraint qualifications or regularity conditions are defined and discussed in \cite{mae1,maciel1}. Through the discussion of this article we consider the Mangasarian-Fromovitz constraint qualification.
\begin{d1}\cite{mfcq1}
The Mangasarian-Fromovitz constraint qualification (MFCQ) is said to be satisfied at a point $x\in\mathbb{R}^n$, if there is a $z\in\mathbb{R}^n$ such that $\nabla g_i(x)^Tz<0$ for $i \in I(x)$.
\end{d1}
Suppose MFCQ holds at $x$. Then the system of inequalities $\nabla g_i(x)^Tz<0$ for $i \in I(x)$ has a nonzero solution $z\in\mathbb{R}^n$. Hence by Gordan's theorem of alternative $\underset{i\in I(x)}{\Sigma} \mu_i \nabla g_i(x)=0$, $\mu_i\geq 0$ has no nonzero solution. That is, $\mu_i=0$ $\forall$ $i\in I(x)$. \\
Conversely suppose the system $\underset{i\in I(x)}{\Sigma} \mu_i \nabla g_i(x)=0,$ $\mu_i\geq 0$ at $x$ has no nonzero solution $\mu$. Then by Gordan's theorem of alternative the system of inequalities $\nabla g_i(x)^Tz<0$ for $i \in I(x)$ has some nonzero solution $z\in\mathbb{R}^n$. \\\\
Above discussion concludes that MFCQ holds at $x$ iff
\bee
\underset{i\in I(x)}{\Sigma} \mu_i \nabla g_i(x)=0, ~~\mu_i\geq 0 ~~\Rightarrow ~~\mu_i=0~~ \forall~~i\in I(x). \lb{MFCQ}
\eee
Strong and weak stationary points for single-objective optimization problems are defined in Definition 1 of \cite{mo1}. In the light of  this definition, strong and weakly critical point of the $MOP$ can be defined, taking care all objective functions together as follows.
\begin{d1}\lb{def_wc}
A feasible point $x$ of the $MOP$ is said to be
\begin{enumerate}[ ( 1)]
\item a strongly critical point of the $MOP$ if  there exist vectors $ \lambda \in \mathds{R}^m_{+}-\{0^m\}$ and $\mu \in \mathds{R}^p_{+}$ satisfying
(\ref{fj1}) and (\ref{fj2}).
\item a weakly critical point of the $MOP$ if  there exists an infeasible sequence $\{x^k\}$ converging to $x\in X$ such that
\bee
\underset{k\rightarrow \infty}{\lim}\frac{\underset{d\in D(x^k)}{\max}~\underset{i\in \Lambda_p}{\max}\{g_i(x^k)+\nabla g_i(x^k)^Td;0\}}{\Phi(x^k)}=1,\lb{wc0}
\eee
where $D(x^k)=\{d: \nabla f_j(x^k)^Td\leq 0,~j\in\Lambda_m\}.$
\end{enumerate}
\end{d1}
One may observe that a strongly critical point is a KKT point of the $MOP$.
\section{A SQP based line search method for MOP}\lb{secsqp}
In order to obtain a feasible descent direction at $x$, we solve a  quadratic programming sub-problem $QP(x)$ at $x$ as
\bea
QP(x):~~\underset{t,d}{\min}~~ t+\frac{1}{2} d^Td ~~& & \nn\\
   \mbox{s. t. }\nabla f_j(x)^Td&\leq& t~~~~j\in\Lambda_m \lb{qcc1}\\
   g_i(x)+\nabla g_i(x)^T d&\leq& t~~~~i\in\Lambda_p.\lb{qcc2}
\eea
Motivation for this sub-problem comes from \cite{flg2} with modifications to address infeasibility. Notice that $t=\Phi(x)$, $d=0$ is a feasible solution of $QP(x)$. Hence the feasibility of the sub-problem is guaranteed. $QP(x)$ has unique solution since this is a convex quadratic problem. The solutions of $QP(x)$ satisfy MFCQ since the system
 \bee
  -\underset{j\in\Lambda_m}{\Sigma} \gamma_j-\underset{i\in\Lambda_p}{\Sigma} \eta_i&=&0\\
   \underset{j\in\Lambda_m}{\Sigma}\gamma_j \nabla f_j(x)+\underset{i\in\Lambda_p}{\Sigma}\eta_i \nabla g_i(x)&=&0\\
   \gamma_j \geq 0, \eta_i\geq 0& &
   \eee
implies $\gamma_j=0$ for all $j$ and $\eta_i=0$ for all $i$. Hence there exist $\lambda\in\mathbb{R}^m_{+}$, $\mu\in\mathbb{R}^p_{+}$, $(\lambda,\mu)\neq 0^{m+p}$ satisfying the KKT optimality conditions. As a result,
\bea
d+\underset{j \in\Lambda_m}{\Sigma} \lambda_j\nabla f_j(x)+\underset{i\in\Lambda_p}{\Sigma} \mu_i \nabla g_i(x)=0, \lb {kd1}\\
1-\underset{j \in\Lambda_m}{\Sigma} \lambda_j-\underset{i\in\Lambda_p}{\Sigma} \mu_i =0 ,\lb{kd2}\\
\lambda_j\geq0,~~~~\lambda_j (\nabla f_j(x)^Td-t)=0,~~~j\in\Lambda_m\lb{kd3},\\
\mu_i\geq0,~~~~\mu_i(g_i(x)+\nabla g_i(x)^Td-t)=0~~~i\in\Lambda_p\lb{kd4},\\
\nabla f_j(x)^Td -t\leq0 ,~~j\in\Lambda_m \lb{kd5},\\
~~g_i(x)+\nabla g_i(x)^Td-t\leq 0 ~~i\in\Lambda_p\lb{kd6}.
\eea
\begin{l1}\lb{cv0}
Suppose that  $(t,d)$ is the solution of the $QP(x)$.
\begin{enumerate}[(I)]
\item Then
\bea
t\leq \Phi(x)-\frac{1}{2}d^Td\lb{bd}.
\eea
\item If $d=0$ and the MFCQ holds at $x$ then $x$ is a strong critical point of MOP.
\item If $d\neq0$ then $d$ is a descent direction of $\Psi_{j,\sigma}(x)$ at $x$ for $\sigma$ sufficiently large.
\end{enumerate}
\end{l1}
\textbf{Proof:}\\
(I) One can observe that $t=\Phi(x)$, $d=0$ is a feasible solution of $QP(x)$. Hence $t+\frac{1}{2}d^Td \leq \Phi(x)$. This implies
\bee
t\leq \Phi(x)-\frac{1}{2}d^Td.
\eee
(II) Suppose $(t,0)$ is the solution of $QP(x)$. Replacing $d$ by $0$ in (\ref{kd1})-(\ref{kd6}), we get
\bea
\underset{j \in\Lambda_m}{\Sigma} \lambda_j \nabla f_j(x)+\underset{i\in\Lambda_p}{\Sigma} \mu_i \nabla g_i(x)=0 \lb {k1}\\
1-\underset{j \in\Lambda_m}{\Sigma} \lambda_j -\underset{i\in\Lambda_p}{\Sigma} \mu_i =0 \lb{k2}\\
\lambda_j\geq0~~~~\lambda_j t=0,~~~j\in\Lambda_m\lb{k3}\\
\mu_i\geq0~~~\mu_i(g_i(x)-t)=0~~~~i\in\Lambda_p\lb{k4}\\
0\leq t,~~~~g_i(x)\leq t ~~~~i\in\Lambda_p\lb{k5}.
\eea
$\Phi(x)\leq t$ follows from definition of $\Phi(x)$ and (\ref{k5}). Then $t$, satisfying (\ref{bd}) with $d= 0$ implies $\Phi(x)\geq t$. Hence $\Phi(x)= t$. From (\ref{k4}), $\mu_i=0$ follows for all $i\notin I(x).$  Also, $\lambda_j>0$ holds for at least one $j$.
Otherwise, (\ref{k1}) and (\ref{k2}) will imply $ \underset{i\in\Lambda_p}{\Sigma} \mu_i \nabla g_i(x)=0$, $\mu_i>0$ for at least one $i$, which violates the MFCQ. This implies that $t=0=\Phi(x)$ (from (\ref{k3})). That is, $x$ is a feasible point. Then from (\ref{k4}),
$\mu_ig_i(x)=0$, $\mu_i\geq0$, $i\in\Lambda_p.$ Therefore, $x$ is a strongly critical point of $P$, which follows from (\ref{k1}).\\
(III) Suppose $(t,d)$ is the solution of $QP(x)$ and $d\neq 0$. Then the following two cases could arise:\\
Case-1: Let $\Phi(x)>0$. Applying (\ref{bd}) in (\ref{qcc2}) we get
\bee
g_i(x)+\nabla g_i(x)^T d\leq t\leq \Phi(x)-\frac{1}{2} d^Td&<&\Phi(x).
\eee
Since $0<\Phi(x)$, we have $\underset{i\in I(x)}{\max}\{g_i(x)+\nabla g_i(x)^T d,0\}-\Phi(x)<0$, from the inequalities above. That is, $\Phi^*(x;d)<0$. If $\sigma$ is chosen in such way that $$\nabla f_j(x)^Td+\sigma \Phi^*(x;d)\leq -\frac{1}{2}{d}^Td<0$$
holds for all $j$ then $d$ will be a descent direction of $\Psi_{j,\sigma}(x)$ for all $j$ (from Lemma 2.1(1) of \cite{baz1}).\\
Case-2: If $\Phi(x)=0$, then $t=0$, $d=0$ is a feasible solution of $QP(x)$. So $g_i(x)+\nabla g_i(x)^Td\leq t\leq 0$ holds for all $i\in I(x)$. So
\bee
\Phi^*(x;d)=\underset{i\in I(x)}{\max}\{g_i(x)+\nabla g_i(x)^Td,0\}-\Phi(x)=0.
 \eee
 Also, $d\neq 0$ implies $t<0$. Hence from (\ref{qcc1}), we have $\nabla f_j(x)^Td\leq t<0$ and consequently $\nabla f_j(x)^Td+\sigma \Phi^*(x;d)<0$.\qed
Let $(t^k,d^k)$ be the solution of the subproblem $QP(x^k)$. Following the arguments of the proof of Lemma \ref{cv0}(III), the penalty parameter $\sigma_k$ can be  updated to force $d^k$ to remain a descent direction for all $\Psi _{j,\sigma_{k+1}}(x^k).$ At $k^{th}$ iteration $\sigma_k$ is unchanged if $d^k$ is the descent direction. Otherwise, $\sigma_k$ is updated as
\bea
\sigma_{k+1}=\max\left\{2\sigma_k,\frac{\nabla f_j(x^k)^Td^k+\frac{1}{2}{d^k}^Td^k}{-\Phi^*(x^k;d^k)},~j\in \Lambda_m \right\}.\lb{sk}
\eea
The theoretical results developed so far are summarized in the following algorithm.
\begin{alg1}\lb{mc1}(A SQP based algorithm)
\begin{enumerate}[Step 1.]
\item (Initialization) Choose $x^{0} \in \mathbb{R}^n$, some scalars $r\in (0,1)$, $\beta\in(0,1)$, the initial penalty parameter $\sigma_0>0$, and an error tolerance $\epsilon$. Set $k:=0$.\lb{initial}
\item Solve the $QP(x^k)$ to find the descent direction $(t^k;d^k)$ with Lagrange multipliers $(\lambda^k$, $\mu^k)$. If $\|d^k\|<\epsilon$, then stop, otherwise go to Step \ref{dneq0}.
\item If $\Phi(x^k)=0$ or $\theta_{j,\sigma_k}(x^k;d^k)\leq-\frac{1}{2} {d^k}^Td^k$ for all $j$, let $\sigma_{k+1}=\sigma_{k}$. Otherwise, $\sigma_{k+1}$ is updated using (\ref{sk}).\lb{dneq0}
\item Compute step length $\alpha_k$ as the first number in the sequence $\{1,r,r^2,...\}$ satisfying \lb{step_alpha}
\bea
\Psi_{j,\sigma_{k+1}}(x^k+\alpha_k d^k)-\Psi_{j,\sigma_{k+1}}(x^k)\leq \alpha_k \beta \theta_{j,\sigma_{k+1}}(x^k;d^k)~~~~\forall~~j\in\Lambda_m.\nn\\ \lb{amj1}
\eea
\item Update $x^{k+1}=x^k+\alpha_k d^k$. Set $k:=k+1$ and go to Step 2.
\end{enumerate}
\end{alg1}
\section{Convergence}\lb{secgconv}
In this section the global convergence of Algorithm \ref{mc1} is proved. The following lemma is used to establish that Step \ref{step_alpha} is well-defined. The extension to the multi-objective case justifies the convergence analysis.
\begin{l1}
Suppose $\nabla f_j(x)$ and $\nabla g_i(x)$ are Lipschitz continuous for every $j\in \Lambda_m$ and $i\in \Lambda_p$ with Lipschitz constant $L$ and let $(t^k,d^k)$ be the solution of the $QP(x^k)$ with $d^k\neq0$. Then (\ref{amj1}) holds for $\alpha$ sufficiently small.
\end{l1}
\textbf{Proof:} Since $\nabla f_j(x)$ and $\nabla g_i(x)$ are Lipschitz continuous for every $j\in \Lambda_m$ and $i\in \Lambda_p$, from Lemma 2.1(3) of \cite{baz1}, there exists $L>0$ such that
\bee
\Psi_{j,\sigma_{k+1}}(x^k+\alpha d^k)\leq \Psi_{j,\sigma_{k+1}}(x^k)+\alpha \theta_{j,\sigma_{k+1}}(x^k;d^k)+\frac{1}{2}(1+\sigma_{k+1})L\alpha^2 \|d^k\|^2
\eee
holds for every $\alpha\in[0,1]$. Hence for every $\alpha\in[0,1]$ and $\beta\in(0,1)$ (initialized in Step \ref{initial} of Algorithm \ref{mc1}) we have,
\bea
& & \Psi_{j,\sigma_{k+1}}(x^k+\alpha d^k)-\Psi_{j,\sigma_{k+1}}(x^k)-\beta\alpha \theta_{j,\sigma_{k+1}}(x^k;d^k) \nn\\
& &\leq (1-\beta) \alpha\theta_{j,\sigma_{k+1}}(x^k;d^k)+\frac{1}{2}(1+\sigma_{k+1})L\alpha^2 \|d^k\|^2.\lb{Lip1}
\eea
Since $d^k\neq 0$, from Step \ref{dneq0} of Algorithm \ref{mc1}, $$(1-\beta)\theta_{j,\sigma_{k+1}}(x^k;d^k)\leq -\frac{1}{2}(1-\beta)\|d^k\|^2<0.$$ Hence from (\ref{Lip1}), (\ref{amj1}) holds for every $\alpha>0$ sufficiently small.\qed
\begin{l1}\lb{cv1}
Let $(t^k,d^k)$ be the solution of the sub-problem $QP(x^k)$ and assume that the sequences $\{x^k\}$ and $\{(t^k,d^k)\}$ are bounded. If $x^k\rightarrow x^*$ as $k\rightarrow \infty$, then $\{(t^k,d^k)\}$ converges to $(t^*,d^*)$, where $(t^*,d^*)$ is the solution of $QP(x^*)$.\\
In particular, if $d^k$ converges to $0$ and the MFCQ holds at every $x^k$ then $x^*$ is a strongly critical point of the $MOP$.
\end{l1}
\textbf{Proof:} If possible let $\{x^k\}$ converges to $x^*$ but $\{(t^k,d^k)\}$ does not converge to $(t^*,d^*)$. Since $\{(t^k,d^k)\}$ is bounded, there exists a sub sequence $\{(t^k,d^k)\}_{k\in K}$ converging to $(\bar{t},\bar{d})\neq(t^*,d^*)$. Since $(t^k,d^k)$ is the optimal solution of $QP(x^k)$, there exists $(\lambda^k,\mu^k)$ such that $(t^k,d^k;\lambda^k,\mu^k)$
satisfies the KKT optimality conditions (\ref{kd1})-(\ref{kd6}).
Now (\ref{kd2}) implies $\{\lambda^k\}$ and $\{\mu^k\}$ are bounded. Hence there exists a converging sub sequence of the subsequence $ \{(\lambda^k,\mu^k)\}_{k \in K}.$  Without loss of generality we may assume $ \lambda^k \rightarrow \lambda^*$ and $\mu^k \rightarrow \mu^*$ as $k\rightarrow \infty$ and $k\in K$. Hence in (\ref{kd1})-(\ref{kd6}),
taking limit $k\rightarrow \infty$, $k \in K$, we have
 \bee
\bar{d}+\underset{j \in\Lambda_m}{\Sigma} \lambda_j^{*}\nabla f_j(x^*)+\underset{i\in\Lambda_p}{\Sigma} \mu_i^{*} \nabla g_i(x^*)=0,  \\
1-\underset{j \in\Lambda_m}{\Sigma} \lambda_j^{*}-\underset{i\in\Lambda_p}{\Sigma}  \mu_i^{*} =0 ,\\
\lambda_j^{*}\geq0~~~~\lambda_j^{*} (\nabla f_j(x^*)^T\bar{d}-\bar{t})=0,~~~j\in\Lambda_m,\\
 \mu_i^{*}\geq0~~~~ \mu_i^{*}(g_i(x^*)+\nabla g_i(x^*)^T\bar{d}-\bar{t})=0~~~i\in\Lambda_p,\\
\nabla f_j(x^*)^T\bar{d} -\bar{t}\leq0 ,~~j\in\Lambda_m ,\\
~~g_i(x^*)+\nabla g_i(x^*)^T\bar{d}-\bar{t}\leq 0 ~~i\in\Lambda_p.
\eee
These imply that $(\bar{t},\bar{d};\lambda^*,\mu^*)$ satisfies first order necessary conditions of the convex quadratic programming sub-problem $QP(x^*)$. Hence $(\bar{t},\bar{d})$ is an optimal solution of $QP(x^*)$. This contradicts the fact that $(t^*,d^*)$ is the optimal solution of $QP(x^*)$, since $QP(x^*)$ has unique solution. Hence $\{(t^k,d^k)\}$ converges to $(t^*,d^*)$. \\
In particular, if $d^k$ converges to $0$ and the MFCQ holds at every $x^k$ then replacing $d^*$ by $0$ in (\ref{kd1})-(\ref{kd6}) at $(x^*,t^*,d^*,\lambda^*,\mu^*)$ and proceeding as in Lemma \ref{cv0}(II) it is easy to prove that $x^*$ is a strongly critical point of the $MOP$.\qed
\begin{l1}\lb{cv3}
Suppose that $\sigma_k=\sigma>0$ for $k$ large enough, the sequences $\{x^k\}$ and $\{(t^k,d^k)\}$ are bounded, $\nabla f_j(x)$ and $\nabla g_i(x)$ are Lipschitz continuous for every $j\in \Lambda_m$ and $i\in \Lambda_p$ with Lipschitz constant $L$, and $\{x^k\}_{k\in K}$ is a convergent subsequence. Then $d^k \rightarrow 0$ as $k \rightarrow \infty$ and $k\in K$.
\end{l1}
\textbf{Proof:} Without loss of generality, assume that $\sigma_k=\sigma$ for all $k\in K$. If possible suppose that there exists an infinite subset $K^{'}\subset K$ and a positive constant $\eta$ such that
\bea
\|d^k\|\geq \eta,~~~~~ k\in K^{'}.\lb{ub1}
\eea
First we will show that there exists $\delta>0$ such that $\alpha_k\geq \delta$ holds for every $k$, where $\alpha_k$ is the step length obtained in Step \ref{step_alpha} of Algorithm \ref{mc1}. From this step either $\alpha_k=1$ or $\alpha_k=r^{k_1}$ holds for some $k_1\in\mathbb{N}$. If $\alpha_k=r^{k_1}$ holds then there exists $\hat{j}\in\Lambda_m$ satisfying,
\bee
\Psi_{\hat{j},\sigma}(x^k+r^{k_1-1} d^k)-\Psi_{\hat{j},\sigma}(x^k)>r^{k_1-1}\beta \theta_{\hat{j},\sigma}(x^k;d^k).
\eee
Then from (\ref{Lip1}),
\bee
\frac{1}{2}(1+\sigma)L r^{2(k_1-1)}\|d^k\|^2\geq -r^{k_1-1}(1-\beta) \theta_{\hat{j},\sigma}(x^k;d^k).
\eee
From Step \ref{dneq0} of Algorithm \ref{mc1},
\bee
\frac{1}{2}(1+\sigma)Lr^{k_1-1}\|d^k\|^2&\geq&\frac{1}{2} (1-\beta)\|d^k\|^2\\
\Rightarrow r^{k_1}&\geq&\frac{r(1-\beta)}{(1+\sigma)L}
\eee
%
Choose $\delta=\min\{1,\frac{r(1-\beta)}{(1+\sigma)L}\}$. Then $\alpha_k\geq \delta$ holds for every $k$. Now
\bee
 \Psi_{j,\sigma}(x^{k+1})-\Psi_{j,\sigma}(x^0)&=&\Sigma_{l=0}^{k} \Psi_{j,\sigma}(x^l+\alpha_l d^l)-\Psi_{j,\sigma}(x^l)\nn\\
 &\leq&-\frac{\beta}{2}\Sigma_{l=0}^{k}  \alpha_k\|d^k\|^2\\
 &\leq&-\frac{1}{2}(k+1)\delta\beta \eta^2~~\forall k\in K'.
\eee
The second inequality follows from (\ref{amj1}) and Step \ref{dneq0} of Algorithm \ref{mc1}. This implies $\Psi_{j,\sigma}(x^k+\alpha_k d^k)\rightarrow -\infty$ as $k\rightarrow \infty$  and $  k\in K'$
(since $\alpha_0 \beta \eta^2>0$). This contradicts the assumption that $\{x^k\}$, $\{(t^k,d^k)\}$ are bounded as $\Psi_{j,\sigma}$ is a continuous function. So there does not exist any $K^{'}\subset K$ and $\eta>0$ such that (\ref{ub1}) holds. Hence the lemma follows.\qed
\begin{l1}\lb{lwc2}
If $\sigma_k \rightarrow \infty$ and the sequences $\{x^k\}$, $\{(t^k,d^k)\}$ are bounded then $\underset{k \rightarrow \infty}{\lim} \Phi(x^k)=0.$
\end{l1}
\textbf{Proof:} Proof of this result is similar to the proof of Lemma 7 in \cite{mo1}.\qed
\begin{t1}
Let $\{x^k\}$ be a sequence generated by Algorithm \ref{mc1}, the sequences $\{x^k\}$ and $\{(t^k,d^k)\}$ are bounded, $\nabla f_j(x)$ and $\nabla g_i(x)$ are Lipschitz continuous for every $j\in \Lambda_m$ and $i\in \Lambda_p$ with Lipschitz constant $L$, and the MFCQ is satisfied at every $x^k$. Then any accumulation point of $\{x^k\}$ is a critical point (either weak or strongly critical point) of the $MOP$.
\end{t1}
\textbf{Proof:}\\
$(i)$ Convergence to a strongly critical point:\\
Let $K$ be an infinite index set such that $x^k\rightarrow x^*$ as $k\rightarrow \infty$ and $ k \in K$. Let $\{(t^k,d^k)\}$ be the solution of $QP(x^k)$. If $d^k\rightarrow 0$ as $k\rightarrow \infty$ then by Lemma \ref{cv1}, $x^*$ is a strongly critical point of $P$.\\
$(ii)$ Convergence to a weakly critical point:\\
If there exists a constant $c_0>0$ such that $\|d^k\|\geq c_0$ for large $k\in K$ then from Lemma \ref{cv3}, $\sigma _k\rightarrow \infty$ as $k\rightarrow \infty$. Since the sequences $\{x^k\}$, $\{(t^k,d^k)\}$ are bounded, so from Lemma \ref{lwc2}, $\underset{k \rightarrow \infty}{\lim} \Phi(x^k)=0$.  Hence from (\ref{bd}) and (\ref{qcc1}), $$\nabla f_j(x^k)^Td\leq t^k< 0~~\forall~~j\in\Lambda_m.$$ This implies, $d^k\in D(x^k) $ for large $k$.
Assume ad absurdum that there exists a constant $\eta_2>0$ such that for sufficiently large $k$,
$$\underset{d\in D(x^k)}{\max}\underset{i \in \Lambda_p}{\max}\{g_i(x^k)+\nabla g_i(x^k)^Td;0\}\leq \Phi(x^k)-\eta_2.$$ Suppose $\hat{d^k}$ maximizes $
\underset{i \in \Lambda_p}{\max}\{g_i(x^k)+\nabla g_i(x^k)^T\hat{d^k},0\}$. Then
\bee
\theta_{j,\sigma_{k+1}}(x^k;d^k)+\frac{1}{2}{d^k}^Td^k&\leq&\Phi(x^k) +\sigma_k\left(\underset{i\in\Lambda_P}{\max}\{g_i(x^k)+\nabla g_i(x^k)^T\hat{d^k},0\}-\Phi(x^k)\right)\\
&\leq&\Phi(x^k)-\eta_2 \sigma_k\\&<&0.
\eee
This contradicts Step 3 of the Algorithm. As a result, $\{x^k\}$ converges to a weakly critical point.
\\
Hence any accumulation point of $\{x^k\}$ is either a strongly critical point or a weakly critical point.\qed
\section{Numerical illustration and discussion}
\begin{figure}[!htbp]
    \centering
    \begin{subfigure}[b]{.49\textwidth}
    \centering
     \includegraphics[height=2.5cm,width=1.1\textwidth]{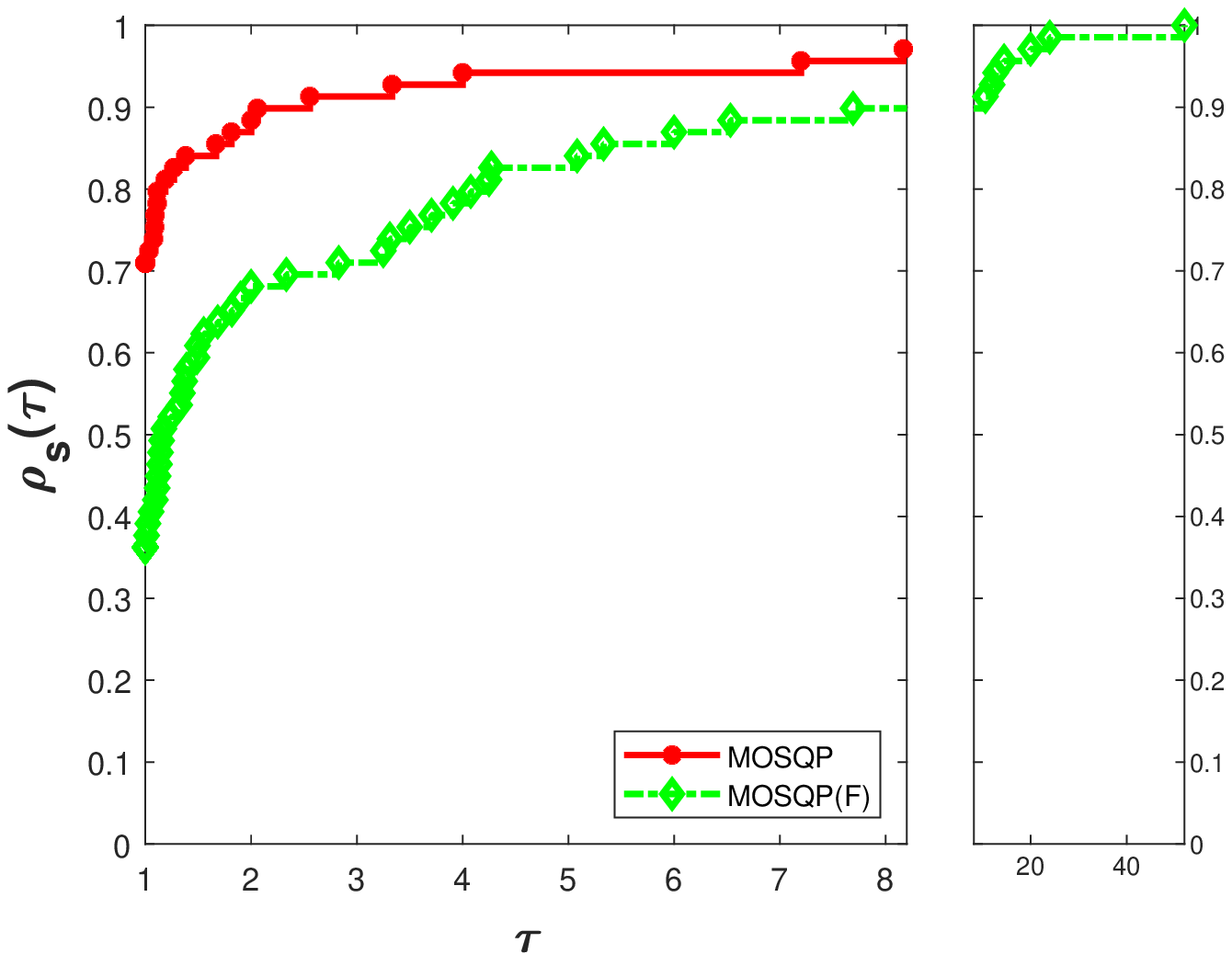}
     \caption{Performance profile between MOSQP and MOSQP(F)}
     \label{pub1}
     \end{subfigure}
     \hfill
    \centering
    \begin{subfigure}[b]{.49\textwidth}
    \centering
     \includegraphics[height=2.5cm,width=1.1\textwidth]{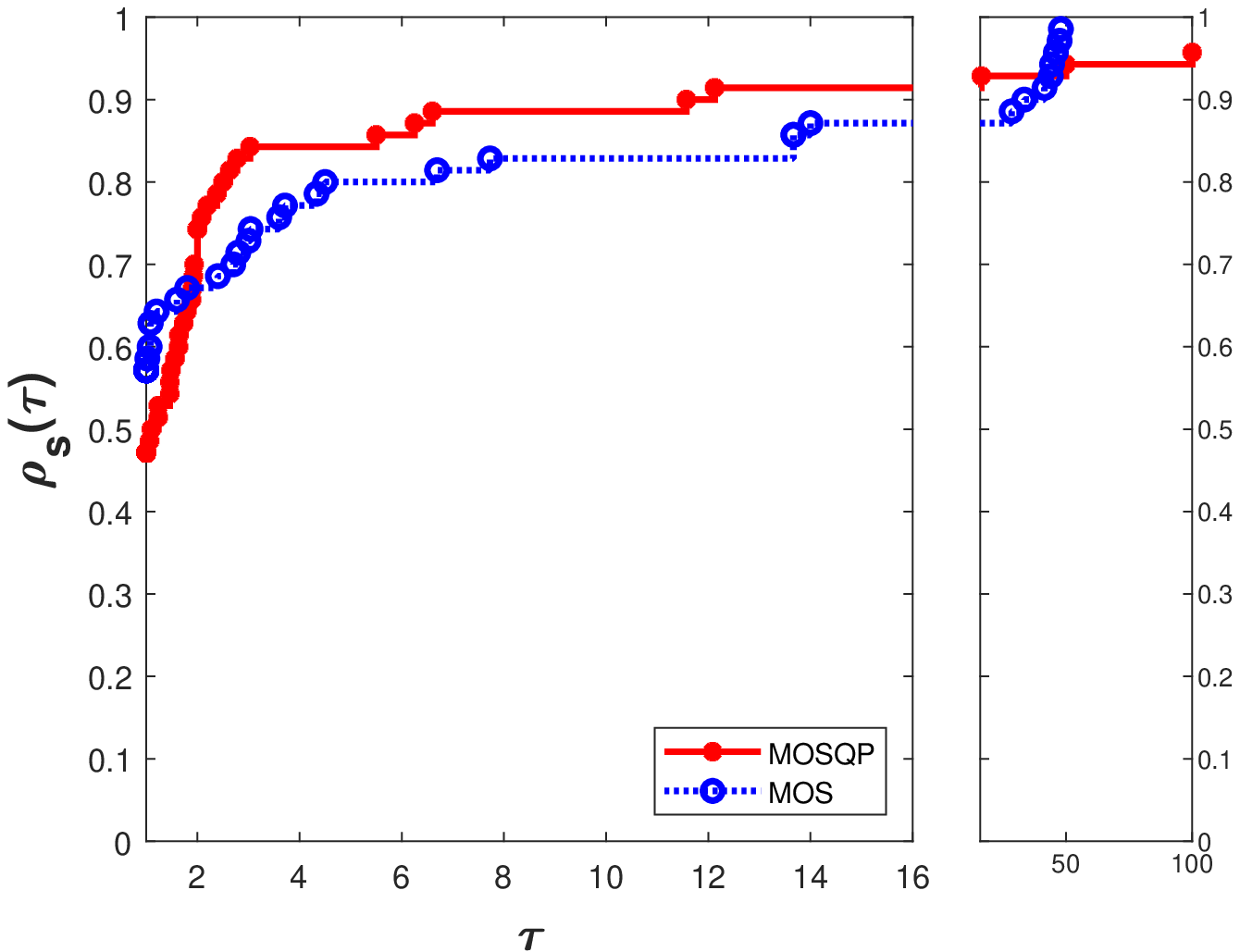}
     \caption{Performance profile between MOSQP and MOS}
     \label{pub2}
     \end{subfigure}
     \caption{Performance profile using purity metric in best run in \textit{\textbf{RAND}}}
\end{figure}

\begin{figure}[!htbp]
    \centering
    \begin{subfigure}[b]{.49\textwidth}
    \centering
     \includegraphics[height=2.5cm,width=1.1\textwidth]{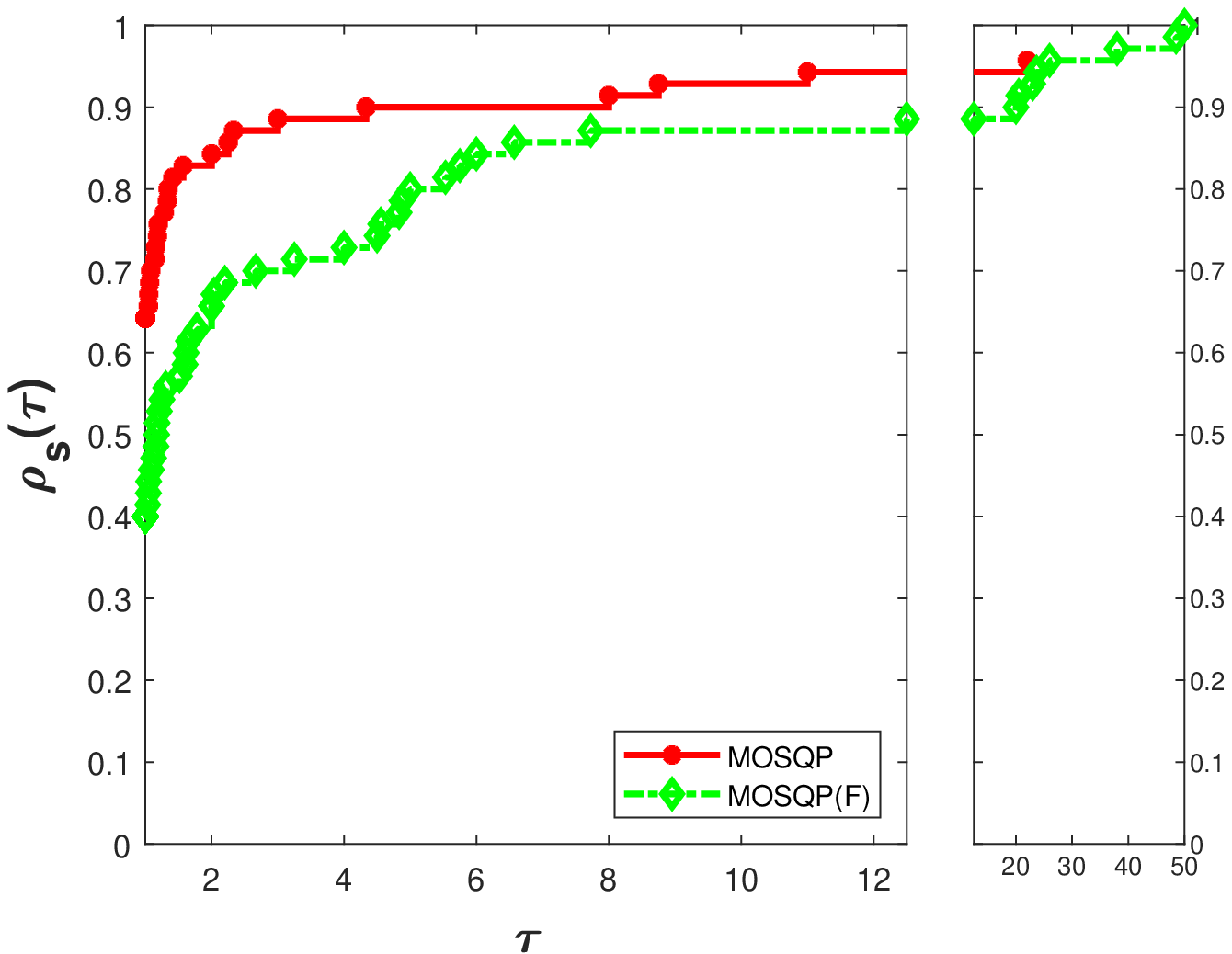}
     \caption{Performance profile between MOSQP and MOSQP(F)}
     \label{puw1}
     \end{subfigure}
     \hfill
    \centering
    \begin{subfigure}[b]{.49\textwidth}
    \centering
     \includegraphics[height=2.5cm,width=1.1\textwidth]{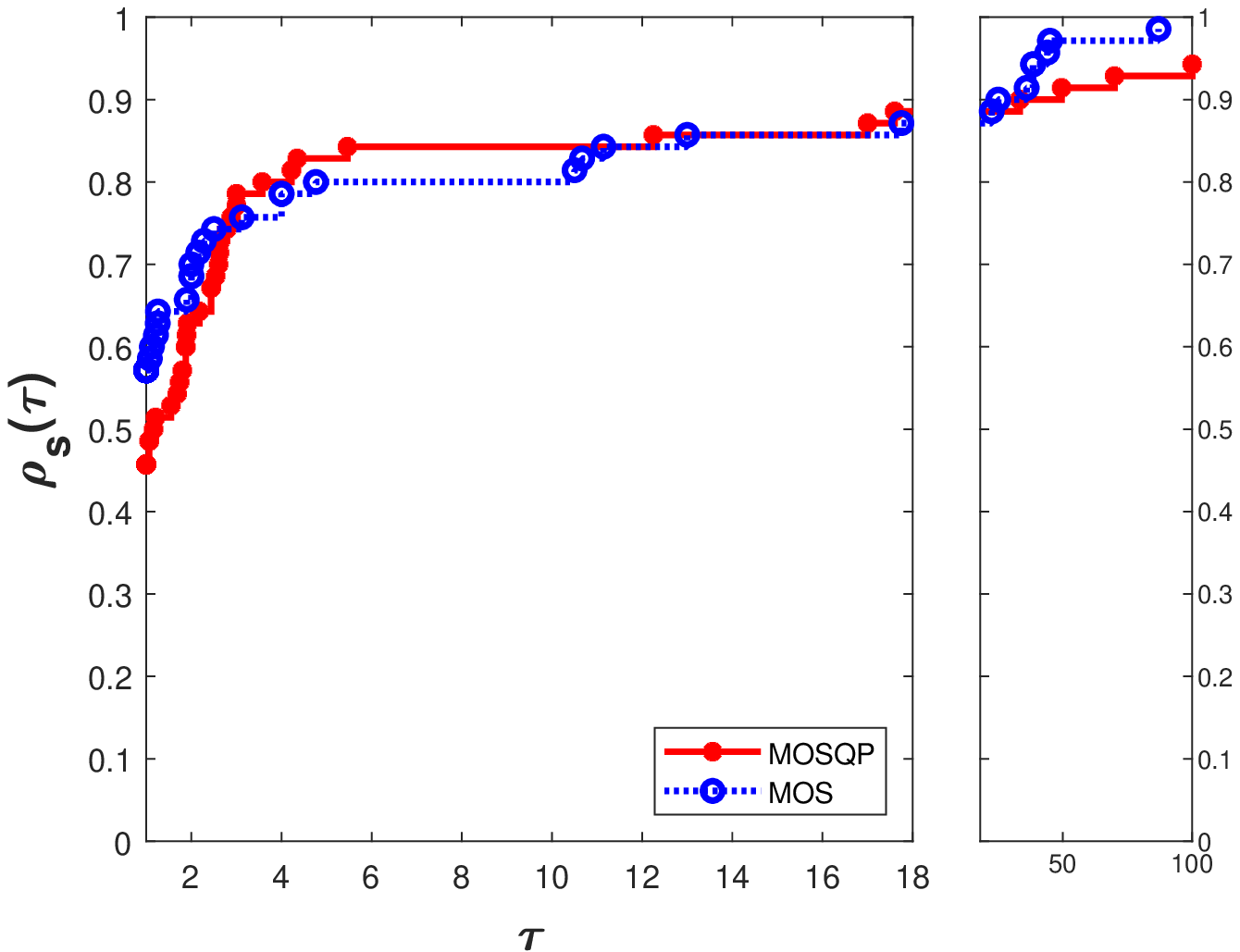}
     \caption{Performance profile between MOSQP and MOS}
     \label{puw2}
     \end{subfigure}
     \caption{Performance profile using purity metric in worst run in \textit{\textbf{RAND}}}
\end{figure}

\begin{figure}[!htbp]
    \centering
    \begin{subfigure}[b]{.49\textwidth}
    \centering
     \includegraphics[height=2.5cm,width=1.1\textwidth]{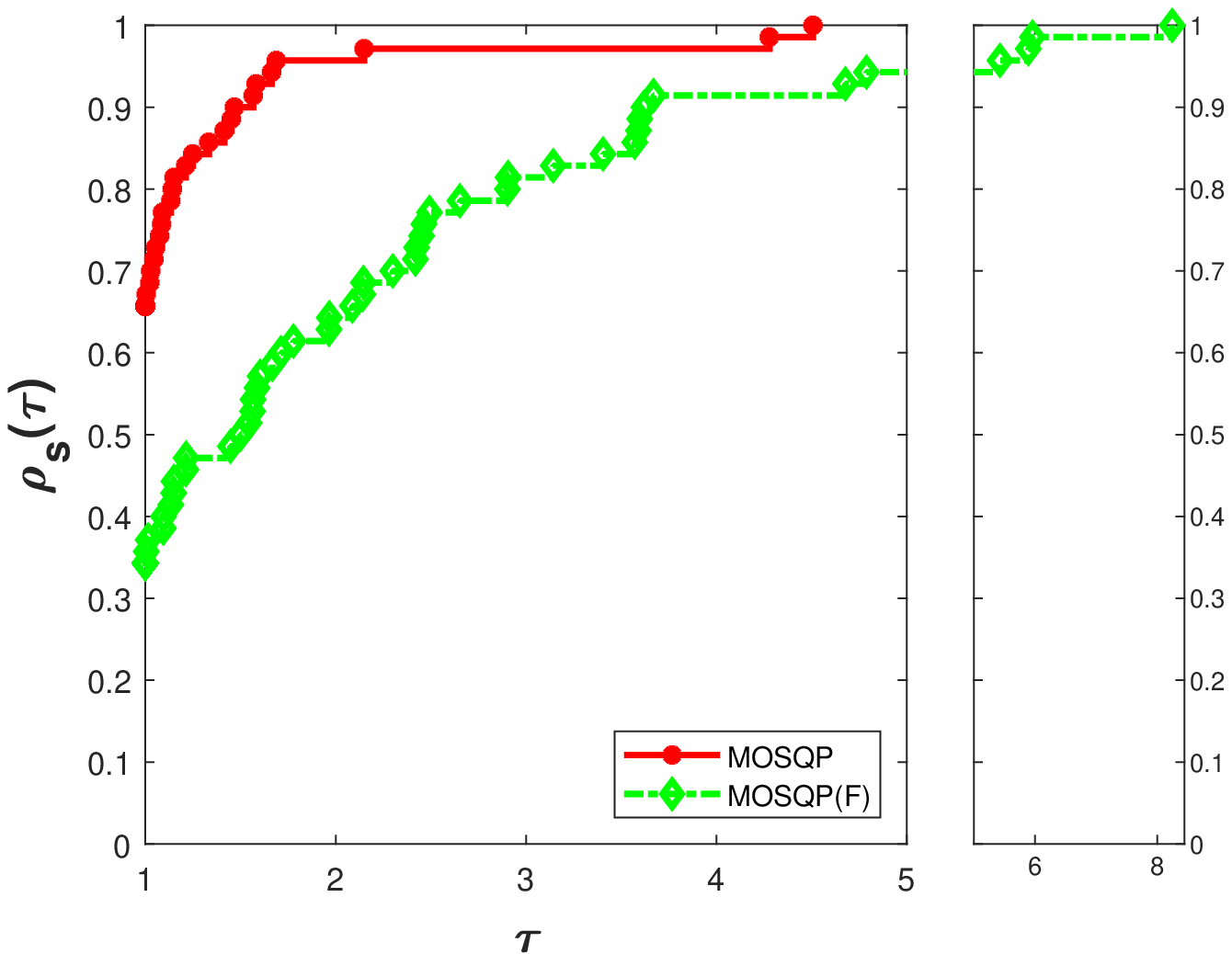}
     \caption{Performance profile between MOSQP and MOSQP(F)}
     \label{gammb1}
     \end{subfigure}
     \hfill
    \centering
    \begin{subfigure}[b]{.49\textwidth}
    \centering
     \includegraphics[height=2.5cm,width=1.1\textwidth]{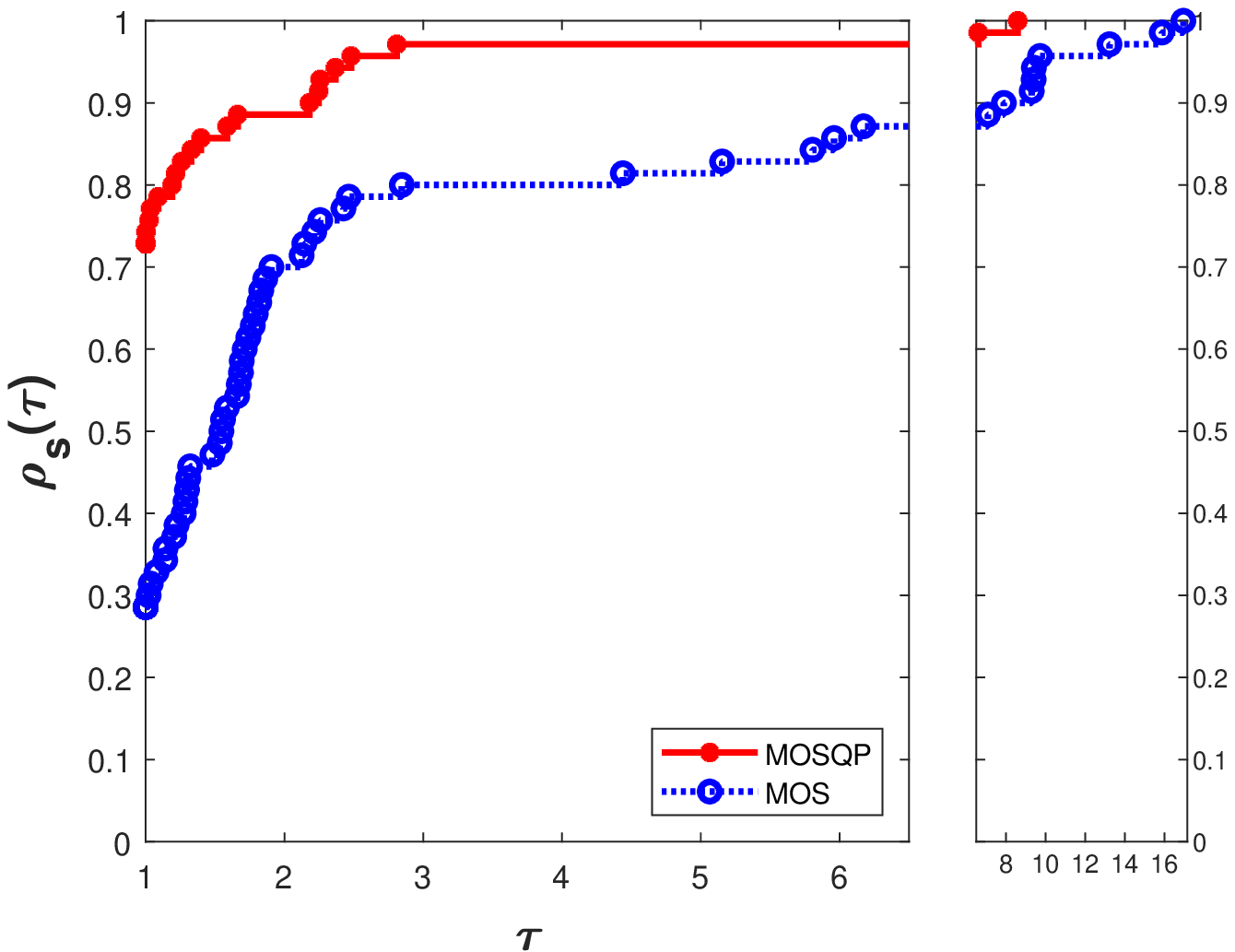}
     \caption{Performance profile between MOSQP and MOS}
     \label{gammb2}
     \end{subfigure}
     \caption{Performance profile using $\Gamma$ metric in best run in \textit{\textbf{RAND}}}
\end{figure}

\begin{figure}[!htbp]
    \centering
    \begin{subfigure}[b]{.49\textwidth}
    \centering
     \includegraphics[height=2.5cm,width=1.1\textwidth]{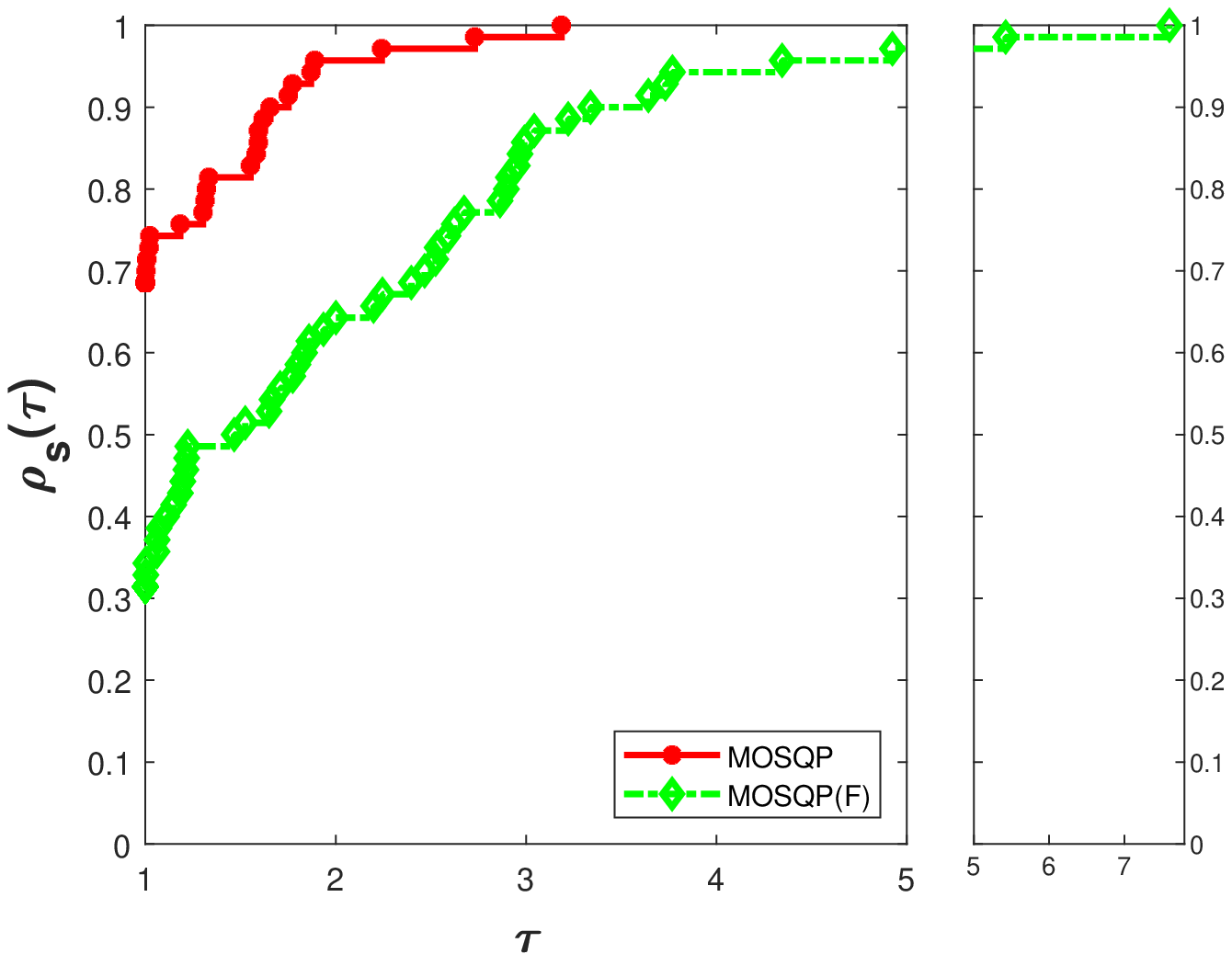}
     \caption{Performance profile between MOSQP and MOSQP(F)}
     \label{gammw1}
     \end{subfigure}
     \hfill
    \centering
    \begin{subfigure}[b]{.49\textwidth}
    \centering
     \includegraphics[height=2.5cm,width=1.1\textwidth]{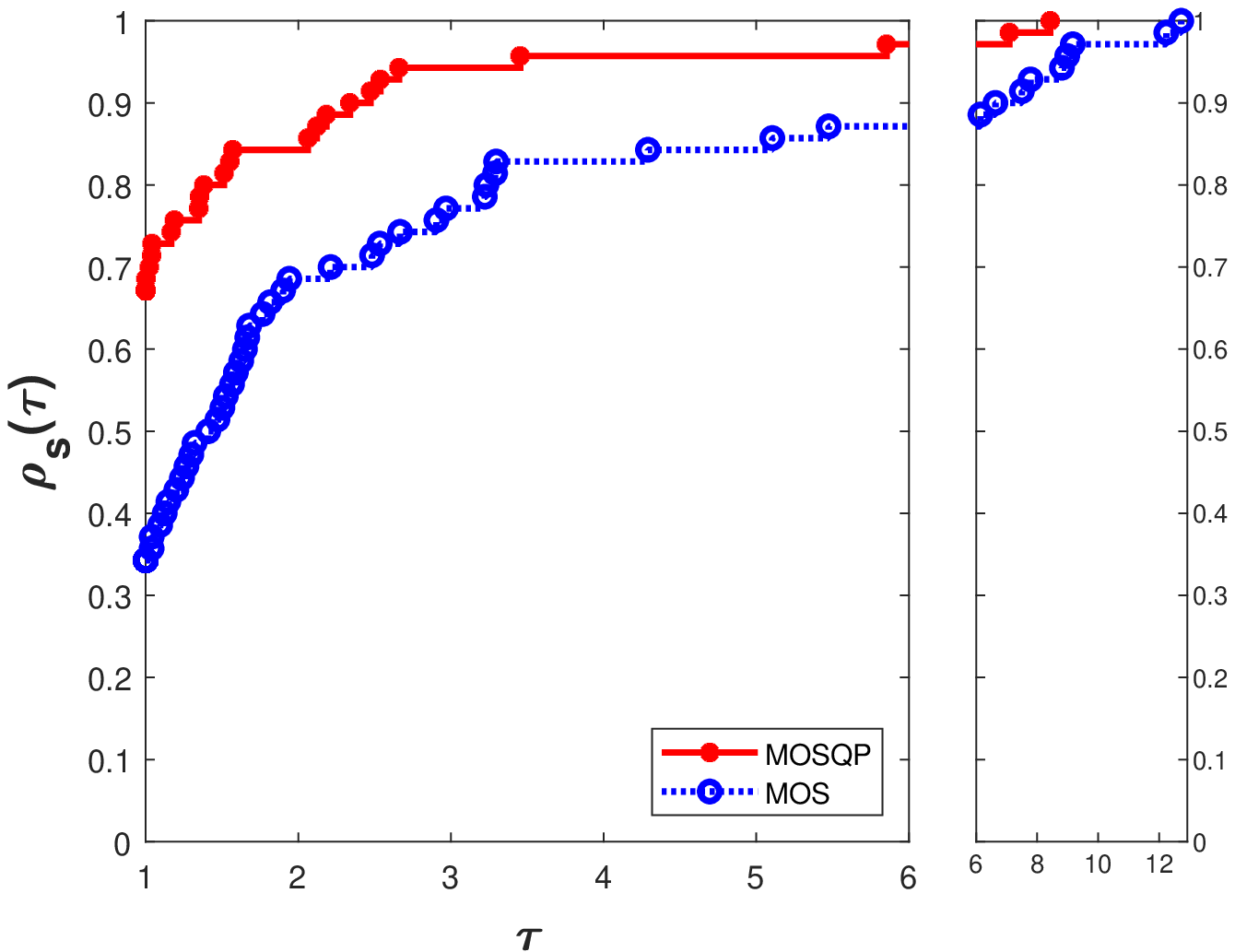}
     \caption{Performance profile between MOSQP and MOS}
     \label{gammw2}
     \end{subfigure}
     \caption{Performance profile using $\Gamma$ metric in worst run in \textit{\textbf{RAND}}}
\end{figure}

\begin{figure}[!htbp]
    \centering
    \begin{subfigure}[b]{.49\textwidth}
    \centering
     \includegraphics[height=2.5cm,width=1.1\textwidth]{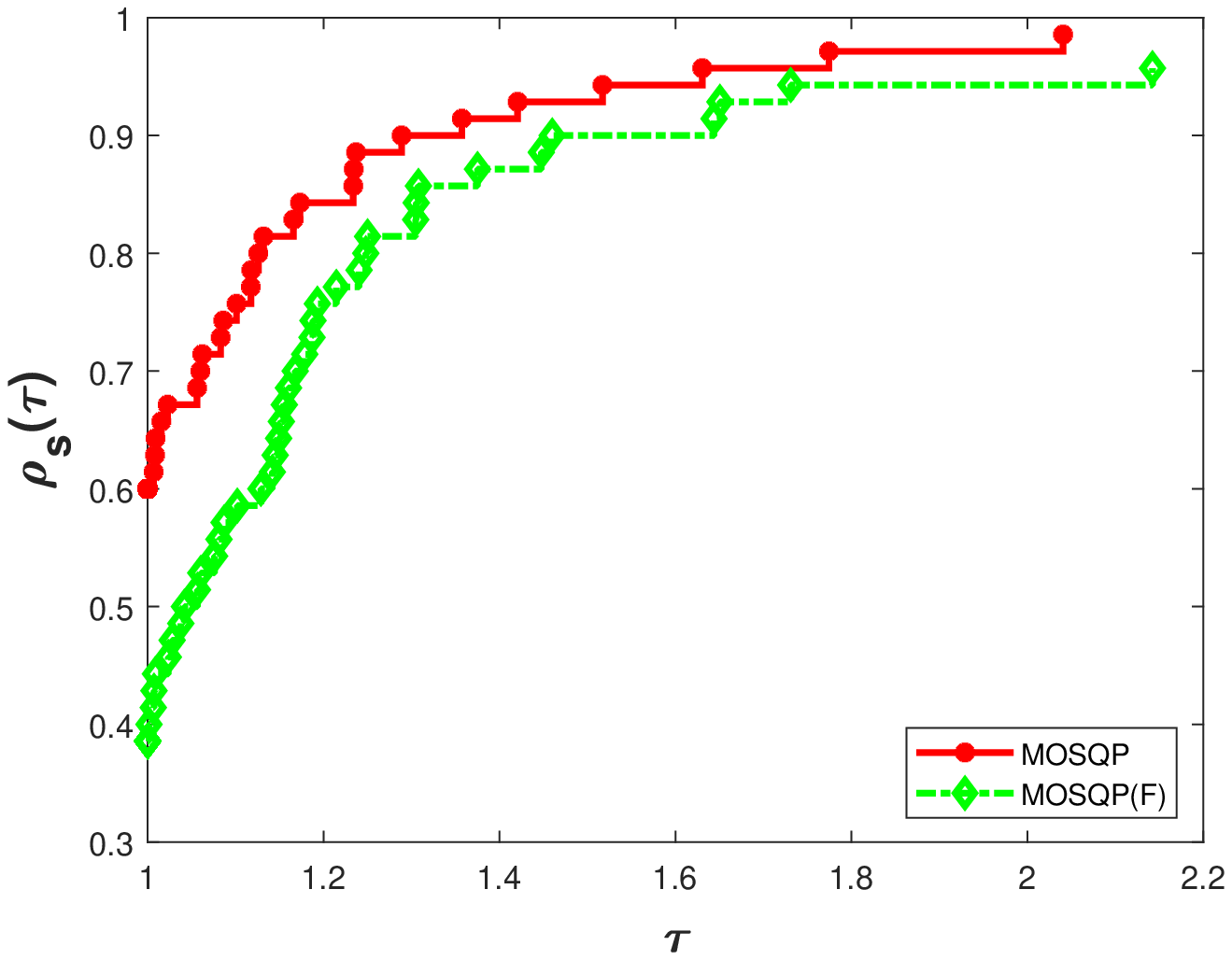}
     \caption{Performance profile between MOSQP and MOSQP(F)}
     \label{delb1}
     \end{subfigure}
     \hfill
    \centering
    \begin{subfigure}[b]{.49\textwidth}
    \centering
     \includegraphics[height=2.5cm,width=1.1\textwidth]{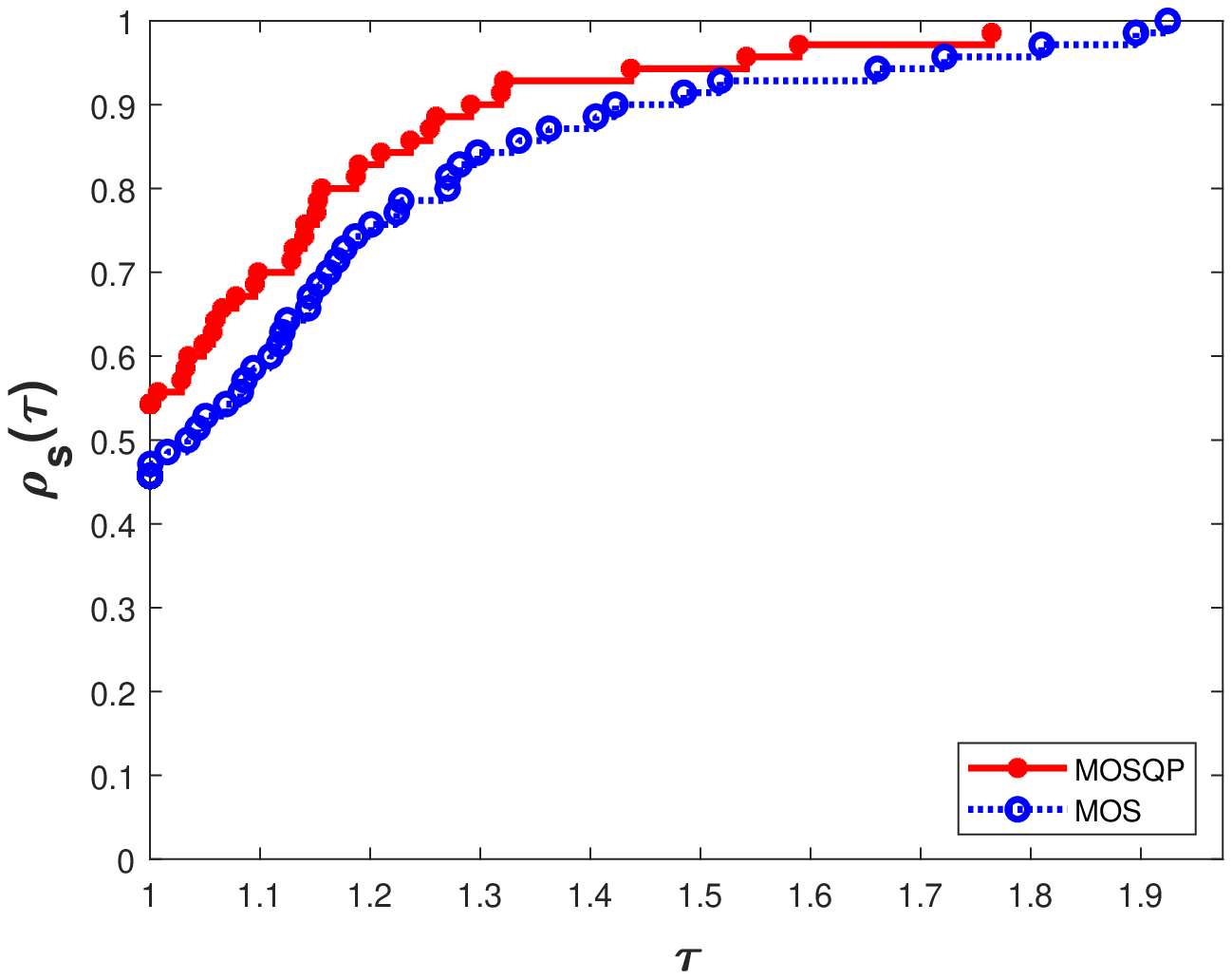}
     \caption{Performance profile between MOSQP and MOS}
     \label{delb2}
     \end{subfigure}
     \caption{Performance profile using $\Delta$ metric in best run in \textit{\textbf{RAND}}}
\end{figure}

\begin{figure}[!htbp]
    \centering
    \begin{subfigure}[b]{.49\textwidth}
    \centering
     \includegraphics[height=2.5cm,width=1.1\textwidth]{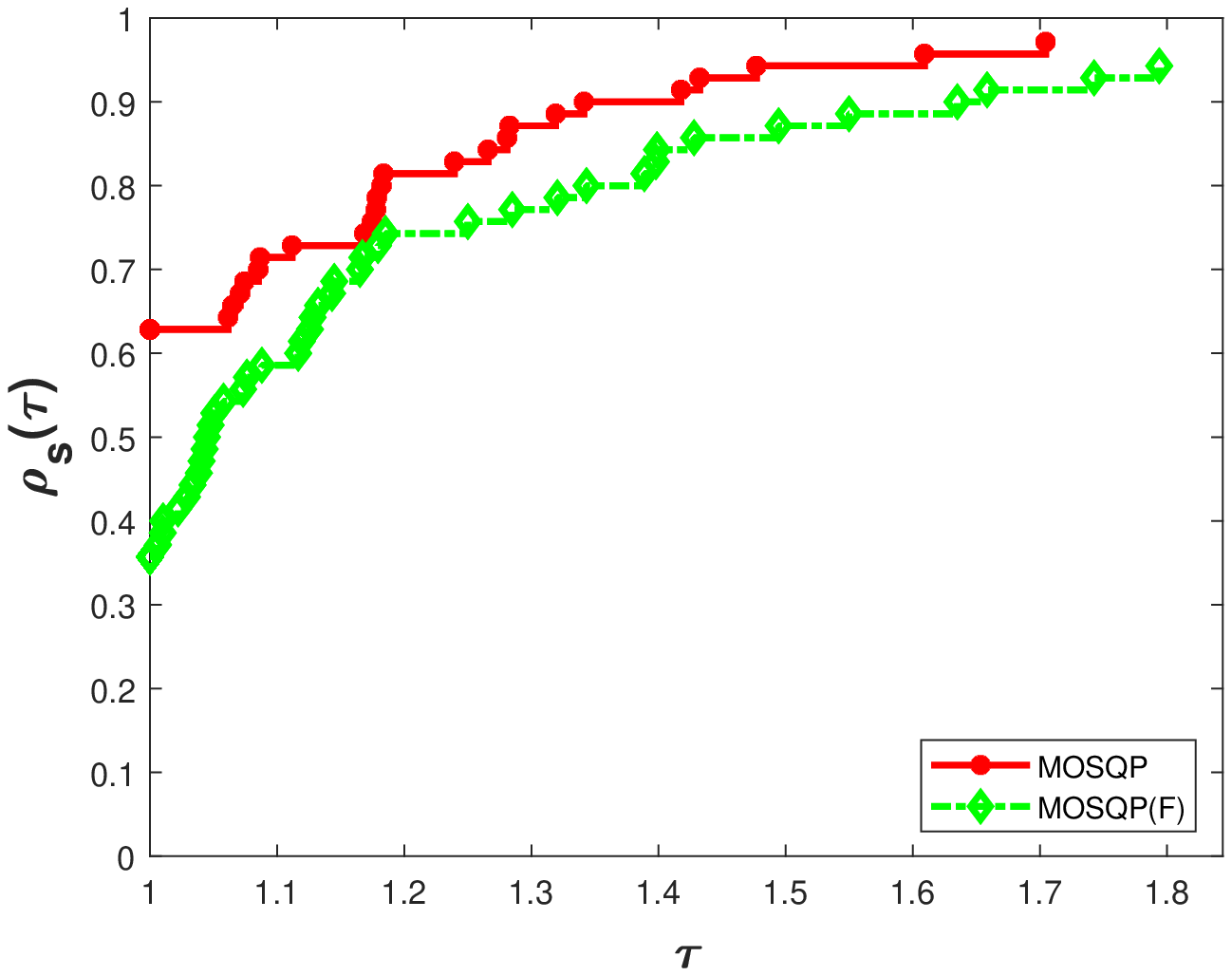}
     \caption{Performance profile between MOSQP and MOSQP(F)}
     \label{delw1}
     \end{subfigure}
     \hfill
    \centering
    \begin{subfigure}[b]{.49\textwidth}
    \centering
     \includegraphics[height=2.5cm,width=1.1\textwidth]{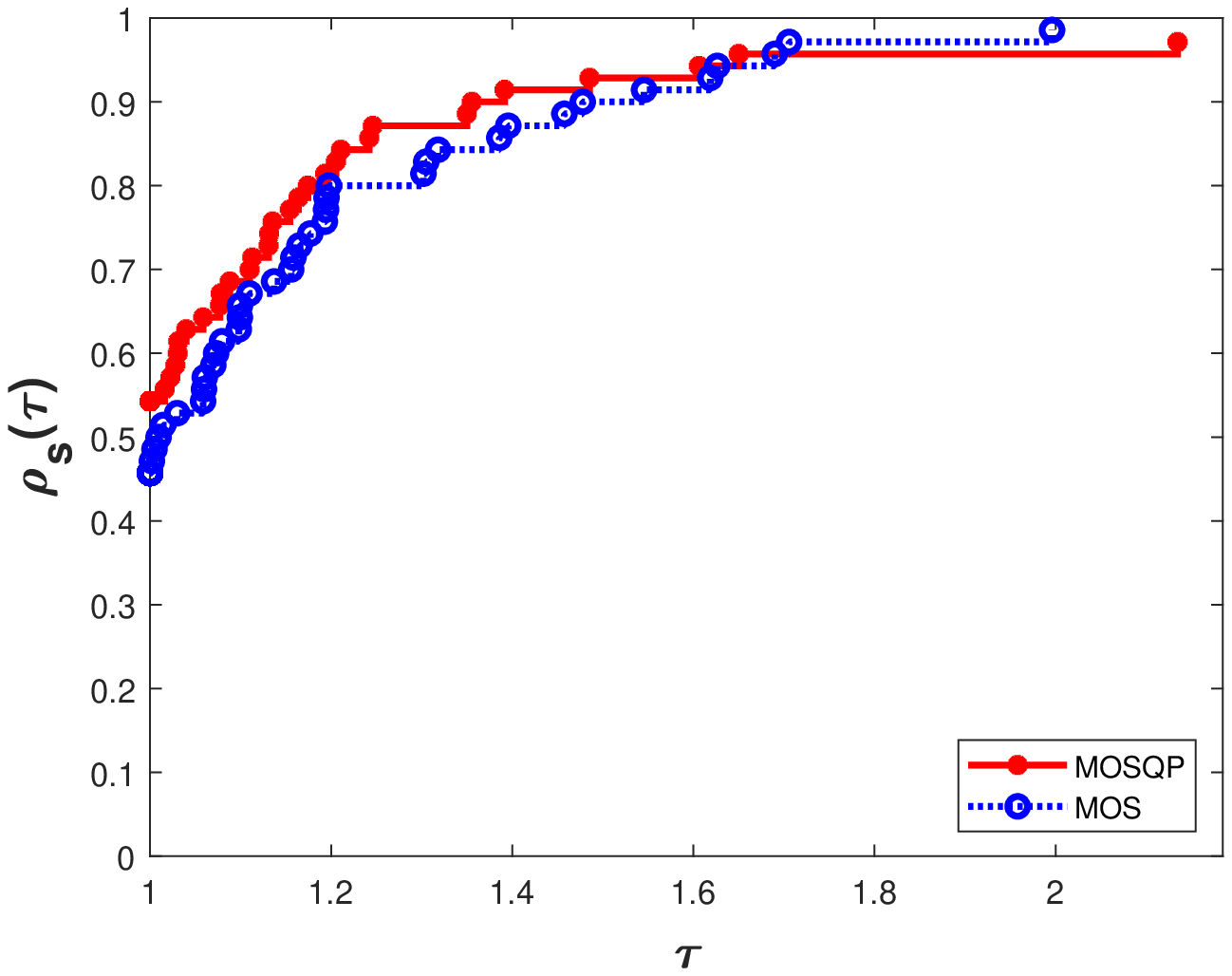}
     \caption{Performance profile between MOSQP and MOS}
     \label{delw2}
     \end{subfigure}
     \caption{Performance profile using $\Delta$ metric in worst run in \textit{\textbf{RAND}}}
\end{figure}

\lb{secnex}
In this section the proposed method (Algorithm \ref{mc1}) (MOSQP) is compared with a classical method (weighted sum method (MOS)) and the method developed by Fliege and Vaz \cite{flg3} (MOSQP(F)). In order to compare different methods we use the performance profiles presented in \cite{flg3,zitq,zitp} with respect to the purity metric and the $\Gamma$ and $\Delta$ spread metrics. (The readers may see the details in \cite{flg3}). In addition to this, two line search techniques MOSQP and MOSQP(F) are compared with respect to average function evaluations.\\\\
\textbf{Performance profile:} Performance profiles are defined by a cumulative function $\rho(\tau)$ representing a performance ratio with respect to a given metric, for a given set of solvers. Given a set of solvers $\mathcal{SO}$ and a set of problems $\mathcal{P}$, let $\varsigma_{p,s}$ be the performance of solver $s$ on solving problem $p$.
The performance ratio is then defined as $r_{p,s}=\frac{\varsigma_{p,s}}{\underset{s\in\mathcal{SO}}{\min}\varsigma_{p,s}}$. The cumulative function $ \rho_s(\tau)$ ($s\in\mathcal{SO}$) is defined as $$ \rho_s(\tau)=\frac{|\{p \in \mathcal{P}:r_{p,s}\leq \tau\}|}{|P|}.$$ It has been observed that the performance profiles are sensitive to the number and types of algorithms considered in the comparison (see \cite{gould1}). So we have compared algorithms pairwise.\\
\textbf{Purity metric:} Purity metric is used to compare the number of non-dominated solutions obtained by different algorithms. Let $F_{p,s}$ be the approximated Pareto front of problem $p$ obtained by method $s$. Then we can build an approximation to the true Pareto front $F_p$ by first considering $\underset {s \in \mathcal{S}}{\cup } F_{p,s}$ and removing the dominated points. The purity metric for algorithm $s$ and problem $p$ is defined by the ratio
$$ \bar{t}_{p,s}=\frac{|F_{p}|}{|F_{p,s}\cap F_{p}|}.$$
Clearly $\bar{t}_{p,s}=\infty$ implies that the algorithm is unable to generate any non-dominated point in the reference Pareto front of the corresponding problem.\\
\textbf{Spread metrics:} Two types of spread metrics ($\Gamma$ and $\Delta$) are used in order to analyze if the points generated by an algorithm are well-distributed in the approximated Pareto front of a given problem. Let $x_1,x_2,...x_N$ be the set of points obtained by a solver $s$ for problem $p$  and let these points be sorted by objective function $j$, i.e., $f_j(x_i)\leq f_j(x_{i+1})$ $(i=1,2,...,N-1)$. Suppose $x_0$ is the best known approximation of global minimum of $f_j$ and $x_{N+1}$ is the best known global maximum of $f_j$, computed over all approximated Pareto fronts obtained by different solvers.
Define $\delta_{i,j}=f_j(x_{i+1})-f_j(x_i)$. Then the $\Gamma$ spread metric is defined by
\bee
\Gamma_{p,s}=\underset{j\in\Lambda_m}{\max}~~\underset{i\in\{0,1,...,N\}}{\max} \delta_{i,j}.
\eee
Define $\bar{\delta_j}$ as the average of the distances $\delta_{i,j}$, $ i=1,2,...,N-1.$ For an algorithm  $s$ and a problem $p$ the spread metric $\Delta_{p,s}$ is
\bee
\Delta_{p,s}=\underset{j\in\Lambda_m}{\max}\left(\frac{\delta_{0,j}+\delta_{N,j}+\Sigma_{i=1}^{N-1} |\delta_{i,j}-\bar{\delta_j}|}{\delta_{0,j}+\delta_{N,j}+(N-1)\bar{\delta_j}}\right).
\eee
\textbf{Test problems:} A set of test problems, collected from different sources, are summarized in Tables \ref{table1} and \ref{table2}. Bound constrained test problems are summarized in Table \ref{table1}. Linear and nonlinear constrained test problems are summarized in Table \ref{table2}. In Table \ref{table2}, `linear' is the number of linear constraints except bound constraints, and `nonlinear' is the number of nonlinear constraints. In both tables $m$ is the number of objective functions and $n$ represents the number of variables.
\begin{table}[!htbp]
\tiny
\begin{center}
\begin{tabular}
{ l c c c|c c c c|c c c  c}\hline
\bf  problem& \bf Source& \bf $m$ & \bf $n$&\bf  problem & \bf Source& \bf $m$ & \bf $n$ &\bf  problem & \bf Source & \bf $m$ & \bf $n$ \\  \hline
BK1&    \cite{hub1}& 2&2  	&Fonseca&\cite{fonseca1}& 2&2 &  MLF1 & \cite{hub1}&2&1 \\ \hline
CEC09\_1&  \cite{zhang1} & 2 &30  &GE2&\cite{adsc1}& 2 &40  & MLF2 & \cite{hub1}&2&2\\ \hline %
CEC09\_2&  \cite{zhang1} & 2 &15  &GE5&\cite{adsc1}& 2 &40   & MOP1& \cite{hub1}& 2&1\\ \hline
CEC09\_3&  \cite{zhang1} & 2 &30  & IKK1&\cite{hub1}& 3&3 & MOP2&   \cite{hub1}& 2&2\\ \hline
CEC09\_7&  \cite{zhang1} & 2 &30  & IM1& \cite{hub1}&2&2  & MOP3&   \cite{hub1}& 2&2\\ \hline
CL1&  \cite{cheng1} & 2 &4 &Jin1& \cite{jin01}& 2&2  & MOP5&   \cite{hub1}& 3&2\\ \hline
Deb41&   \cite{deb1999multi}& 2&2  &Jin2\_a&   \cite{jin01}&  2&2& MOP6&  \cite{hub1}&  2&2  \\ \hline
Deb513&   \cite{deb1999multi}& 2&2  & Jin3&  \cite{jin01}&  2&2 & MOP7&  \cite{hub1}&  3&2  \\ \hline
Deb521a\_a&   \cite{deb1999multi}&2&2& Jin4\_a&   \cite{jin01}&  2&2 & SK1&    \cite{hub1}&  2&1 \\ \hline
Deb521b&   \cite{deb1999multi}&  2&2   & KW2&\cite{adsc1}& 2&2 & SK2&  \cite{hub1}&  2&4\\ \hline
DG01&    \cite{hub1}& 2&1   &lovison1& \cite{lovi1}& 2&2 & SP1&   \cite{hub1}& 2&2    \\ \hline
DTLZ1&   \cite{debs1}& 3&7  & lovison2&\cite{lovi1}&  2&2  & SSFYY1&    \cite{hub1}& 2&2\\ \hline
DTLZ1n2&   \cite{debs1}& 2&2  & lovison3&\cite{lovi1}&  2&2  & SSFYY2&  \cite{hub1}& 2& 1\\ \hline
DTLZ2&   \cite{debs1}& 3&12  & lovison4&\cite{lovi1}&  2&2  & TKLY1&    \cite{hub1}& 2&4\\ \hline
DTLZ2n2&   \cite{debs1}& 2&2  &lovison5& \cite{lovi1}& 3&3 & VFM1&  \cite{hub1}& 3& 2 \\ \hline
DTLZ5\_a&   \cite{debs1}& 3&12  & lovison6& \cite{lovi1}& 3&3 &  VU1&    \cite{hub1}&2&2 \\ \hline
DTLZ5n2\_a&   \cite{debs1}& 2&2  &LRS1&  \cite{hub1}&2&2    &   VU2&    \cite{hub1}&2&2 \\ \hline
ex005&   \cite{hwang2}& 2&2 & MHHM1&    \cite{hub1}& 3&1  &ZDT3&    \cite{zitzler1}&2&30 \\ \hline
Far1& \cite{hub1}& 2&2   &   MHHM2&     \cite{hub1}& 3&2    &  &   & &  \\ \hline
\end{tabular}
\caption { Multi-objective test problems with bound constraints}
\lb{table1}
\end{center}
\end{table}
\begin{table}[!htbp]
\tiny
\begin{center}
\begin{tabular}
{l c c c c c| c c c c c c}\hline
\bf  problem&\bf Source & \bf $m$ & \bf $n$ & \bf linear& \bf nonlinear&
\bf  problem& \bf Source & \bf $m$ & \bf $n$ & \bf linear& \bf nonlinear \\  \hline
ABC\_Comp&  \cite{hwang1}&2  &2&2&1   & GE3& \cite{adsc1}& 2&2&0&2 \\ \hline
BNH&   \cite{deb1999multi} &2& 2&0&2   & GE4&\cite{adsc1}& 3&3&0&1 \\ \hline
CEC09\_C3&    \cite{zhang1} &2& 10&0&1   &liswetm &\cite{leyffer1} &2&7&5&0 \\ \hline
CEC09\_C9&    \cite{zhang1} &3& 10&0&1   &MOQP\_002& \cite{leyffer1}&3&20&9&0 \\ \hline
ex003 & \cite{tappeta}& 2 &2 &0 & 2  &OSY &  \cite{deb1999multi}& 2&6&4&2   \\ \hline
ex004&  \cite{oliveira1}&  2&2&2&0 & SRN &  \cite{deb1999multi} & 2&2&1&1  \\ \hline
GE1& \cite{adsc1}& 2&2&0&1&TNK & \cite{deb1999multi}&  2&2&0&2 \\ \hline
\end{tabular}
\caption { Multi-objective test problems with linear and nonlinear constraints }
\lb{table2}
\end{center}
\end{table}
\\\\{\bf Implementation details:} MATLAB code (2019a) is developed  for all three methods. The MATLAB code of MOSQP(F) is available in public domain, which is not used here. For MOSQP(F), we have developed own code which uses only the Step 4 (third stage) of Algorithm 4.1 \cite{flg3} since the convergence analysis of Algorithm 4.1 \cite{flg3} is different from the convergence analysis of MOSQP. Multi start techniques, similar to MOSQP, is used to generate an approximated Pareto front for MOSQP(F).
\begin{itemize}
\item Quadratic sub-problems are solved using MATLAB function `\textit{quadprog}' with `\textit{Algorithm}',`\textit{interior-point-convex}'.
\item For MOS, the test problems are converted to single objective optimization problems and solved using MATLAB function `\textit{fmincon}' with '\textit{Algorithm}'  `\textit{sqp}', Specified `\textit{objective gradient}' and `\textit{constraint gradient}', and initial approximation as $(l+u)/2$, where $l$ and $u$ are used as in \cite{flg3}.
\item $\|d^k\|<10^{-5}$ or a maximum of $500$ iterations are considered as stopping criteria.
\item It is essential to find a set of well distributed solutions of $MOP$. Spreading out an approximation to a Pareto set is a difficult problem. One simple technique may not work always in a satisfactory manner for all type of problems. Here, to generate an approximated Pareto front, we have selected the initial point with the strategies \textit{\textbf{LINE}} and \textit{\textbf{RAND}} and random parameters in the scalarization method. \textit{\textbf{LINE}} is considered only for bi-objective optimization problems and \textit{\textbf{RAND}} is considered for both bi-objective and more than two objective optimization problems.
\begin{itemize}
\item Initial point selection strategy \textit{\textbf{LINE}} is considered for bi-objective optimization problems. Here $100$ initial points are chosen in the line segment joining $l$ and $u$, i.e. $x^{0,k}=l+k\frac{u-l}{99}$, $k=0,1,2,...,99$ and for MOS we have solved problems of the form\\$ \underset{x\in X} {\min}~~ w f_1(x)+(1-w)f_2(x)$ for $w=k/99$, $k=0,1,2,...,99$.
\item For every test (two or three objective) problem initial points selection strategy \textit{\textbf{RAND}} is considered. Here $100$ random initial points are selected uniformly distributed  in $l$ and $u$, and for MOS  we have solved\\$\underset{x\in X} {\min}~~  \underset{j\in\Lambda_m}{\Sigma} w_j f_j(x)$ $w_j\geq 0$ $w\neq 0$, where $w$ is a random vector. Every test problem is executed $10$ times with random initial points and weights.
\end{itemize}
  \item Restoration procedure is not used for MOSQP(F) if the quadratic sub-problem is infeasible. These points are excluded. Quadratic sub-problem ($QP(x^k)$)  of Algorithm \ref{mc1} always has a solution since this is a convex quadratic problem and has at least one feasible solution.
\item Different run with initial point selection strategy \textit{\textbf{RAND}} generates different set of non-dominated points. Among 10 runs the run which generates highest number of non-dominated solutions, is denoted as best run. Similarly, the run which generates lowest number of non-dominated solutions, is denoted as worst run. Performance profiles are compared for best and worst runs.
\end{itemize}
The performance profiles between MOSQP and MOSQP(F) using purity metric of best run in \textit{\textbf{RAND}} is provided in Figure \ref{pub1} and the performance profiles between MOSQP and MOS using purity metric of best run in \textit{\textbf{RAND}} is provided in Figure \ref{pub2}. Figures \ref{puw1} and \ref{puw2} correspond to the performance profiles for the purity metric in worst run comparing MOSQP with MOSQP(F) and MOSQP with MOS, respectively. The performance profiles for the $\Gamma$ metric in best run comparing MOSQP with MOSQP(F) and MOSQP with MOS are provided in Figures \ref{gammb1} and \ref{gammb2} respectively. The performance profiles for the $\Gamma$ metric in worst run comparing MOSQP with MOSQP(F) and MOSQP with MOS are provided in Figures \ref{gammw1} and \ref{gammw2} respectively. Figures \ref{delb1} and \ref{delb2} correspond to the performance profiles for $\Delta$ metric in best run comparing MOSQP with MOSQP(F) and MOSQP with MOS, respectively. The performance profiles for the $\Delta$ metric in worst run comparing MOSQP with MOSQP(F) and MOSQP with MOS are provided in Figures \ref{delw1} and \ref{delw2} respectively.\\\\

\begin{figure}[!htbp]
    \centering
    \begin{subfigure}[b]{.49\textwidth}
    \centering
     \includegraphics[height=2.5cm,width=1.1\textwidth]{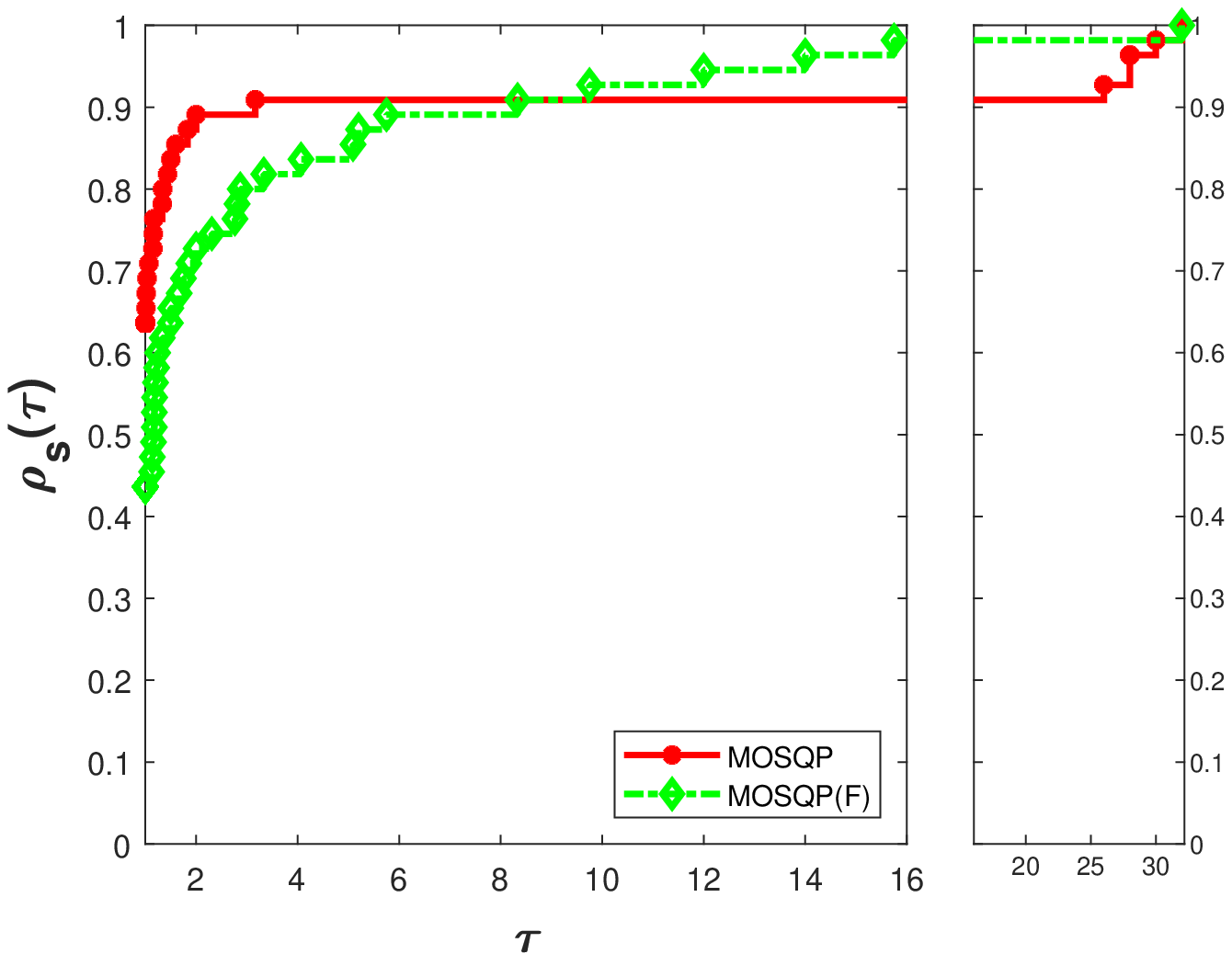}
     \caption{Performance profile between MOSQP and MOSQP(F)}
     \label{pul1}
     \end{subfigure}
     \hfill
    \centering
    \begin{subfigure}[b]{.49\textwidth}
    \centering
     \includegraphics[height=2.5cm,width=1.1\textwidth]{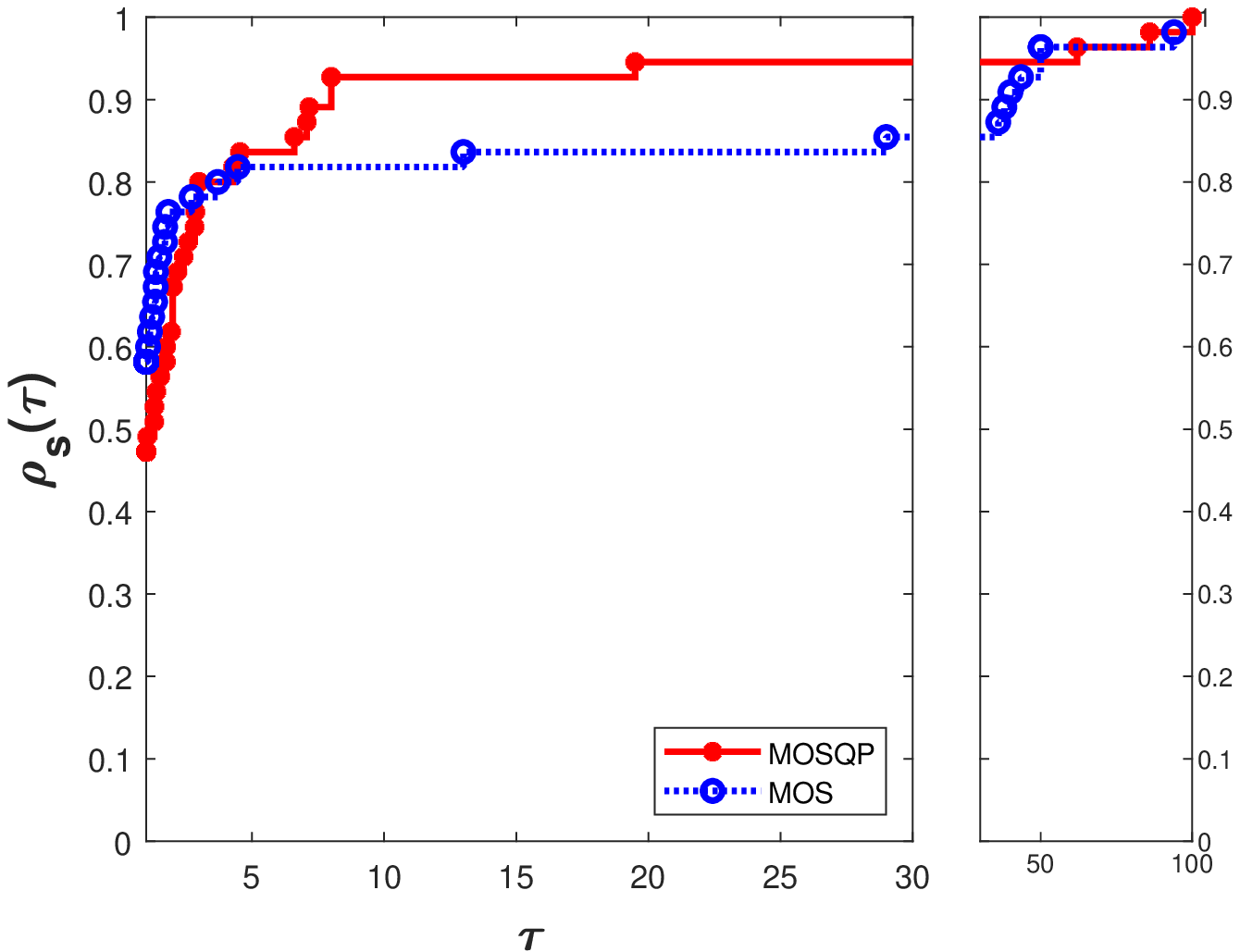}
     \caption{Performance profile between MOSQP and MOS}
     \label{pul2}
     \end{subfigure}
     \caption{Performance profile using purity metric in \textit{\textbf{LINE}}}
\end{figure}

\begin{figure}[!htbp]
    \centering
    \begin{subfigure}[b]{.49\textwidth}
    \centering
     \includegraphics[height=2.5cm,width=1.1\textwidth]{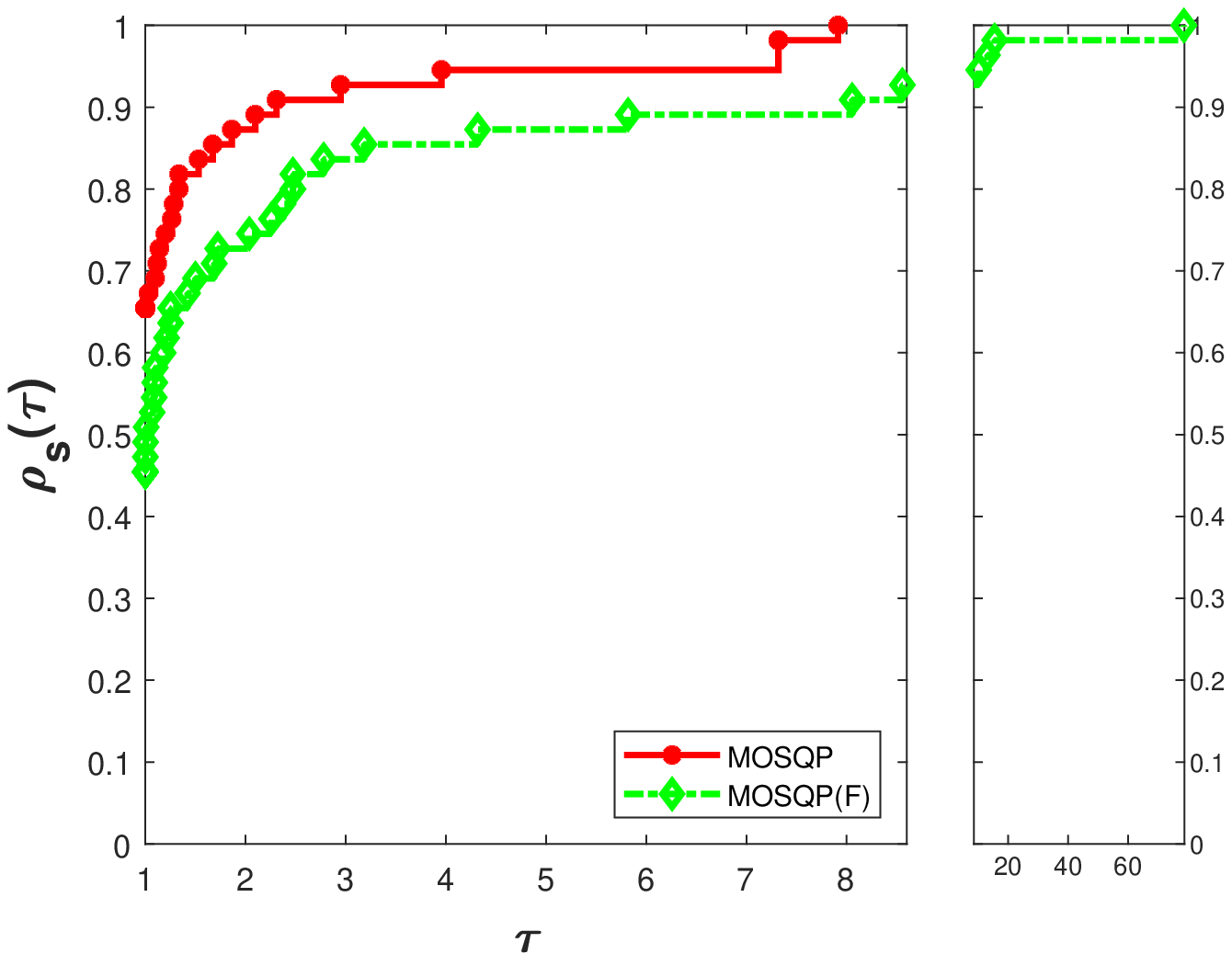}
     \caption{Performance profile between MOSQP and MOSQP(F)}
     \label{gammal1}
     \end{subfigure}
     \hfill
    \centering
    \begin{subfigure}[b]{.49\textwidth}
    \centering
     \includegraphics[height=2.5cm,width=1.1\textwidth]{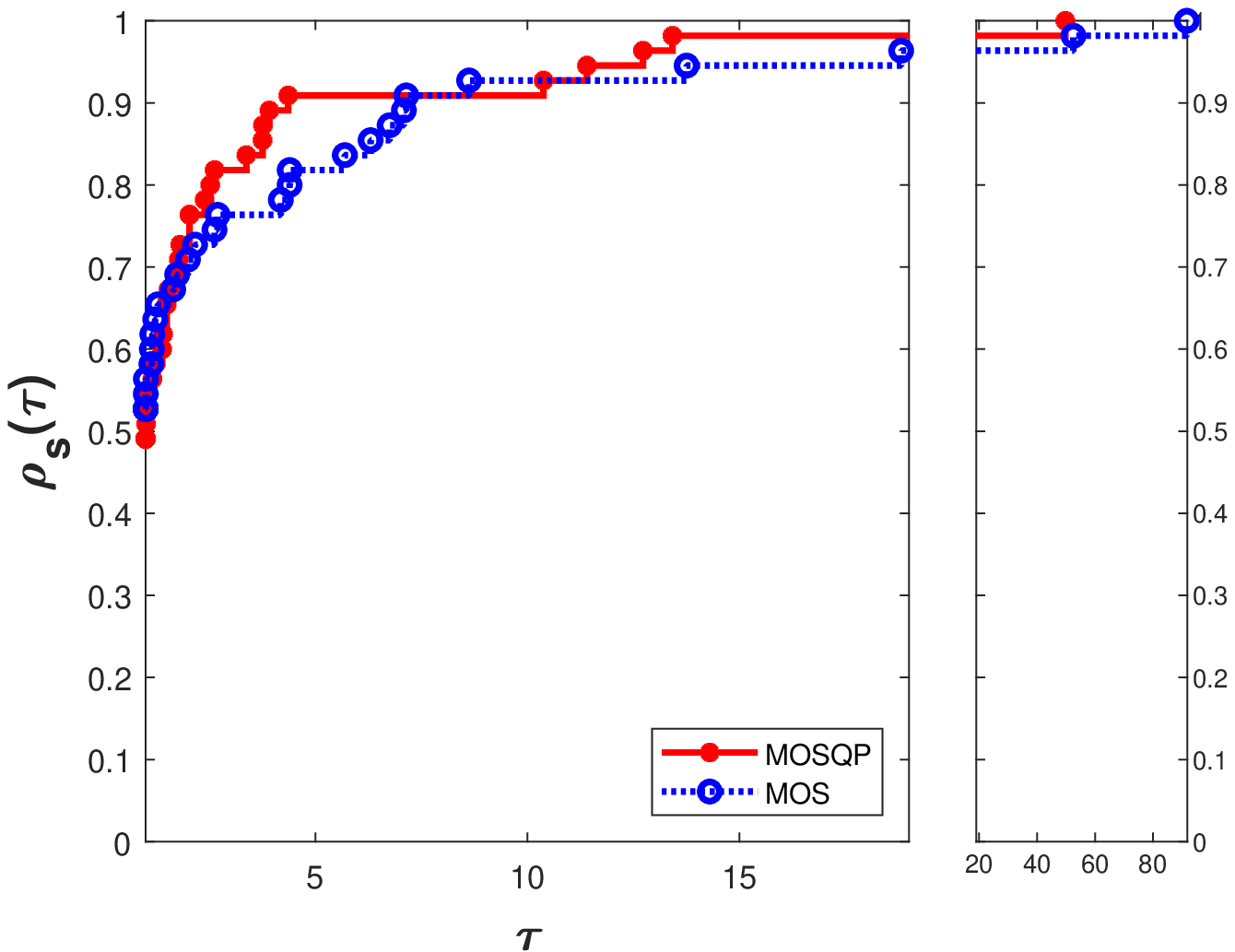}
     \caption{Performance profile between MOSQP and MOS}
     \label{gammal2}
     \end{subfigure}
     \caption{Performance profile using $\Gamma$ metric in \textit{\textbf{LINE}}}
\end{figure}

\begin{figure}[!htbp]
    \centering
    \begin{subfigure}[b]{.49\textwidth}
    \centering
     \includegraphics[height=2.5cm,width=1.1\textwidth]{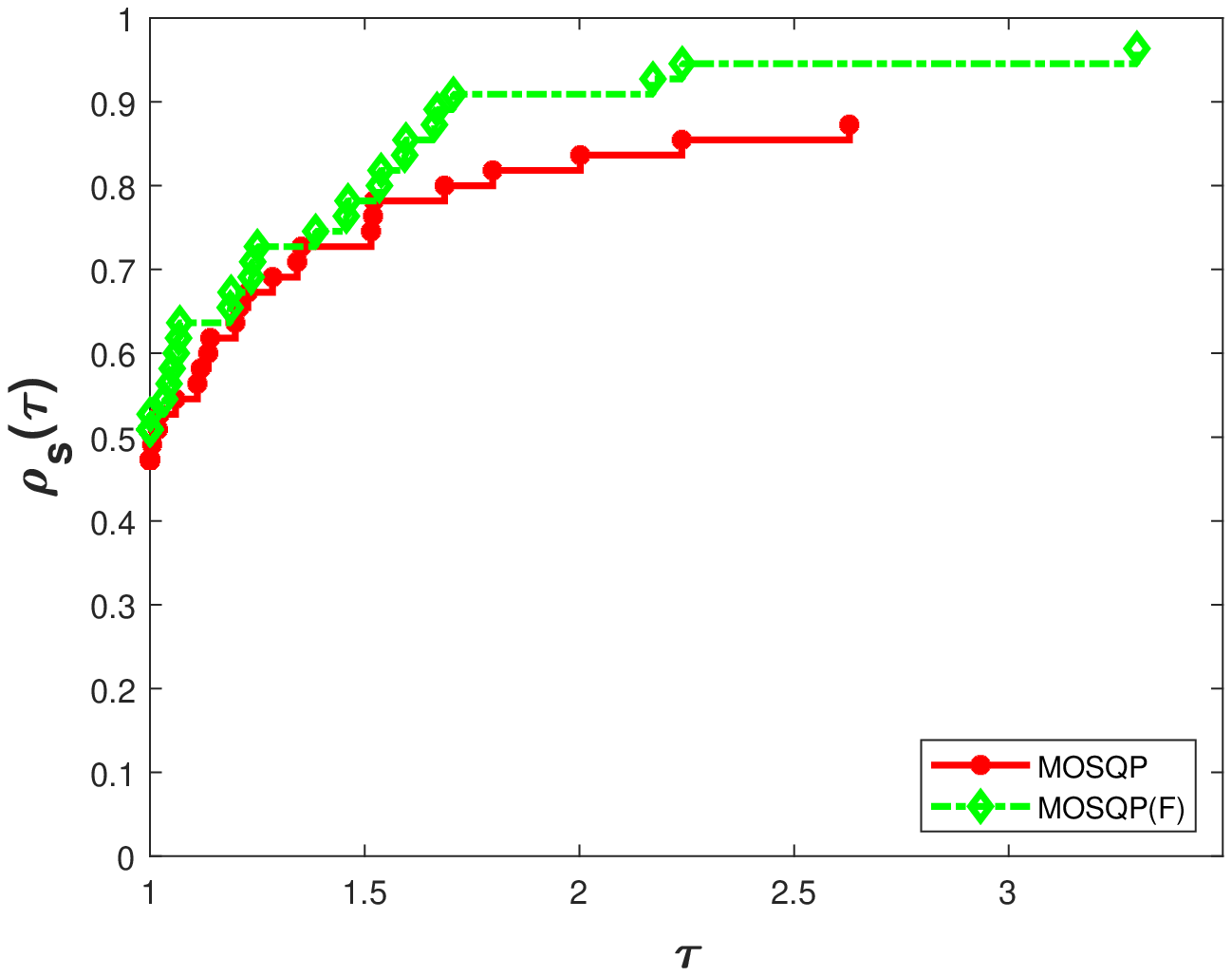}
     \caption{Performance profile between MOSQP and MOSQP(F)}
     \label{dell1}
     \end{subfigure}
     \hfill
    \centering
    \begin{subfigure}[b]{.49\textwidth}
    \centering
     \includegraphics[height=2.5cm,width=1.1\textwidth]{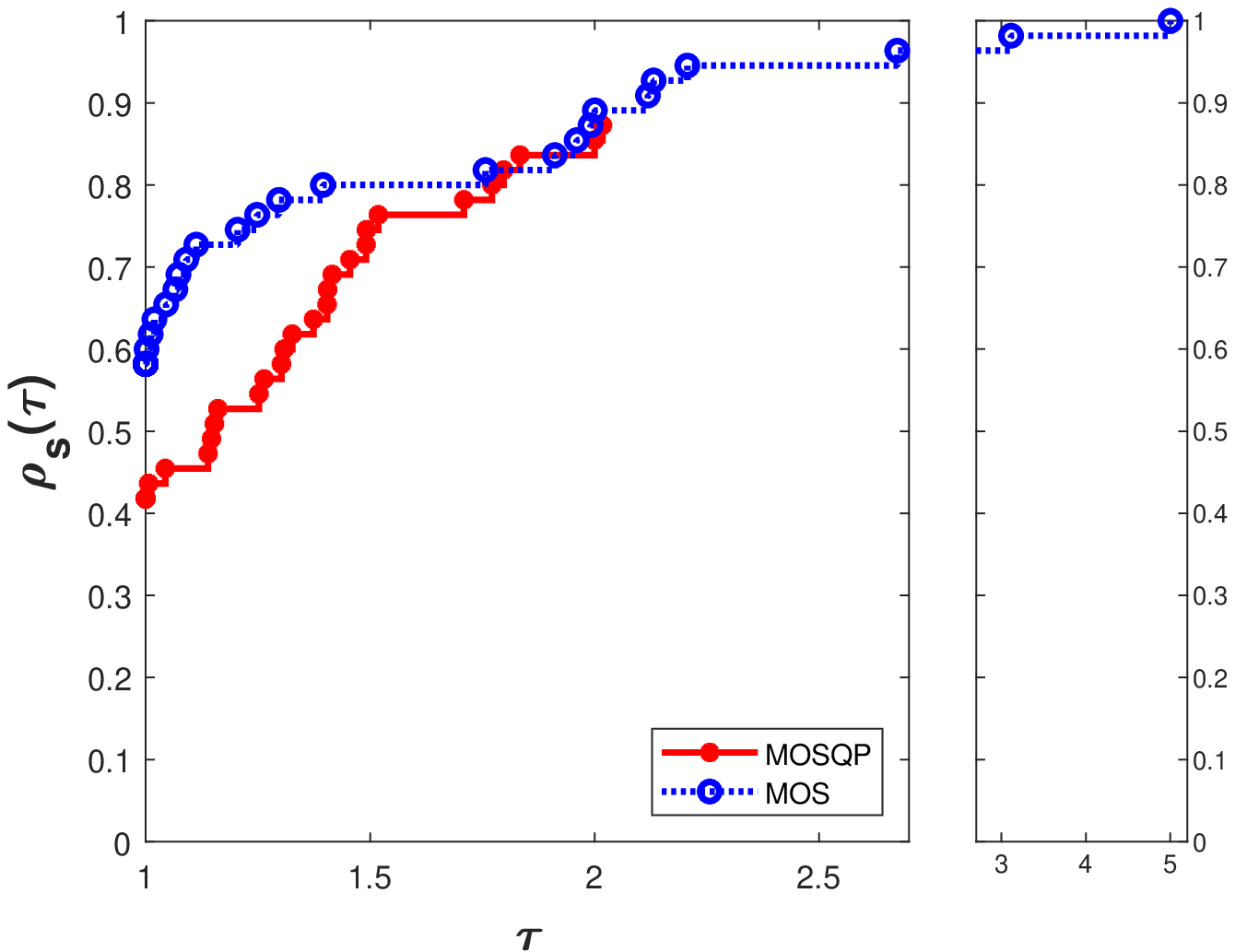}
     \caption{Performance profile between MOSQP and MOS}
     \label{dell2}
     \end{subfigure}
     \caption{Performance profile using $\Delta$ metric in \textit{\textbf{LINE}}}
\end{figure}

\begin{figure}[!htbp]
    \centering
    \begin{subfigure}[b]{.49\textwidth}
    \centering
     \includegraphics[height=2.5cm,width=1.1\textwidth]{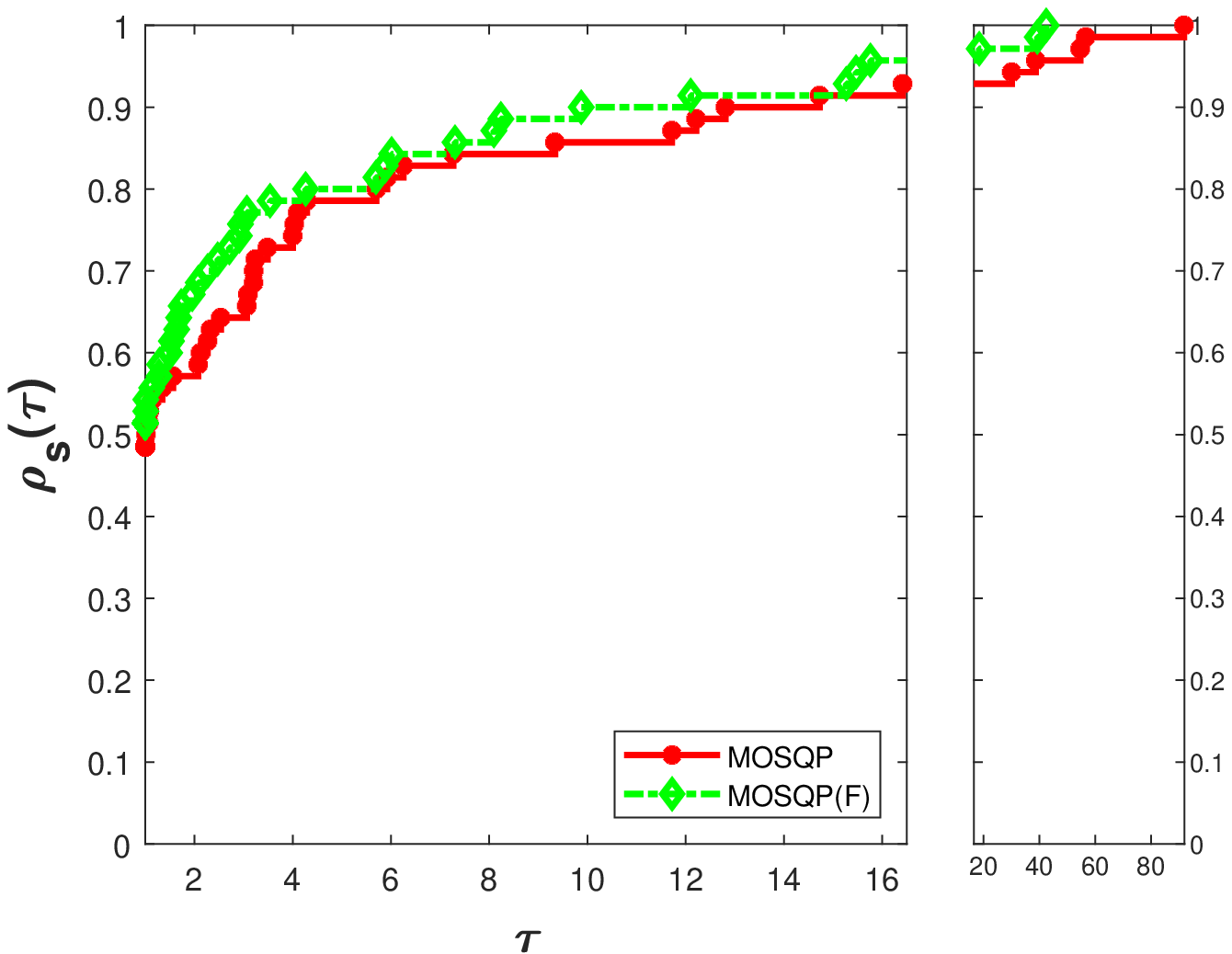}
     \caption{Performance profile for best run}
     \label{ctb}
     \end{subfigure}
     \hfill
    \centering
    \begin{subfigure}[b]{.49\textwidth}
    \centering
     \includegraphics[height=2.5cm,width=1.1\textwidth]{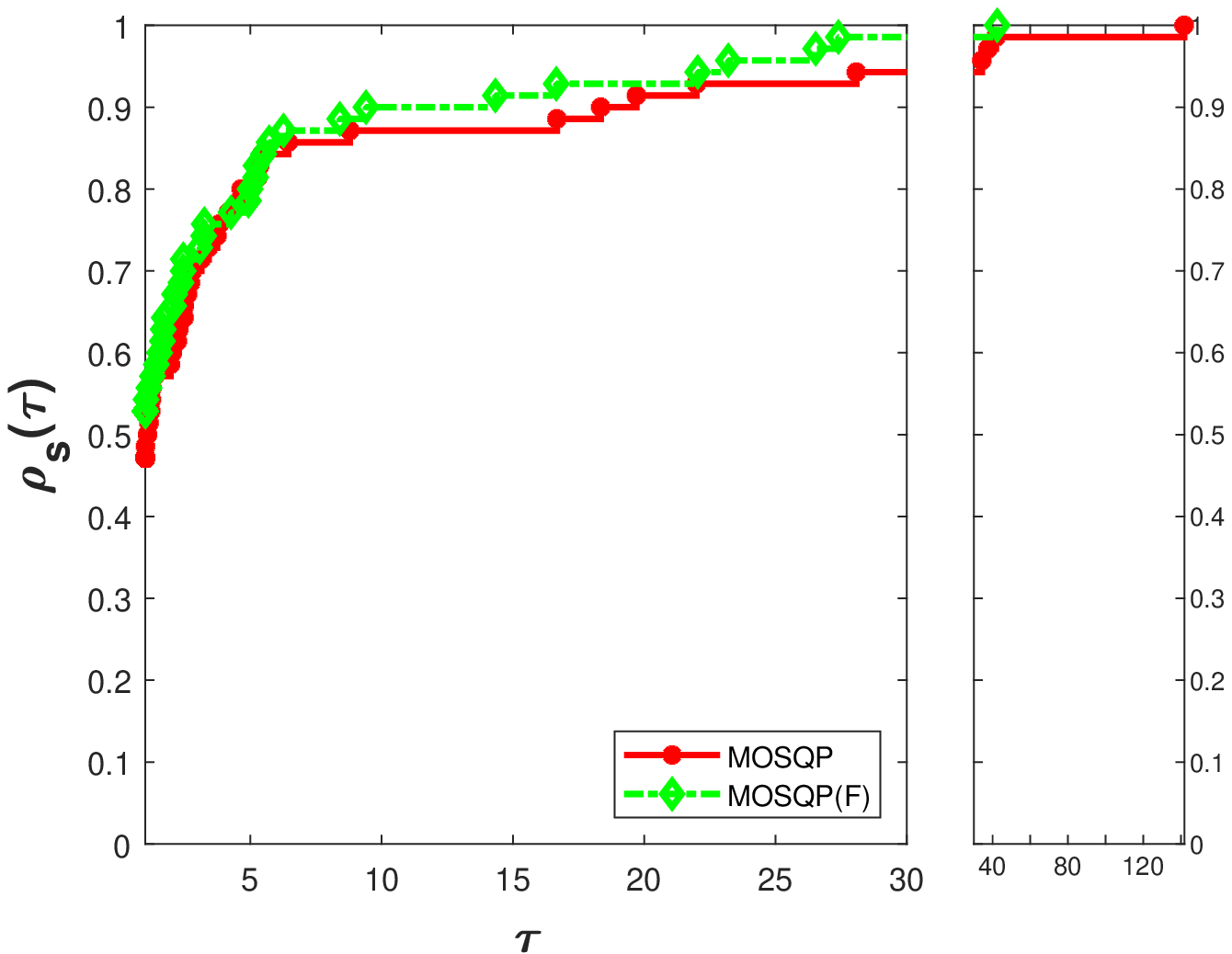}
     \caption{Performance profile for worst run}
     \label{ctw}
     \end{subfigure}
     \caption{Performance profile between MOSQP and MOSQP(F) using average function evaluations in \textit{\textbf{RAND}}}
\end{figure}

\newpage
The performance profiles using purity metric in \textit{\textbf{LINE}} between MOSQP and MOSQP(F) and between MOSQP and MOS are provided in Figures \ref{pul1} and \ref{pul2} respectively. Figures \ref{gammal1} and \ref{gammal2} correspond to the performance profiles for $\Gamma$ metric in \textit{\textbf{LINE}} comparing MOSQP with MOSQP(F) and MOSQP with MOS, respectively. The performance profiles for the $\Delta$ metric in \textit{\textbf{LINE}} comparing MOSQP with MOSQP(F) and MOSQP with MOS are provided in Figures \ref{dell1} and \ref{dell2} respectively.\\\\
\\
Two line search techniques MOSQP and MOSQP(F) are compared using average number of function evaluations per non-dominated points. We have calculated gradient for MOSQP and MOSQP(F) and Hessian for MOSQP(F) using forward difference formula. If $n_1$ and $n_2$ are number of non-dominated points generated by MOSQP and MOSQP(F), then average function evaluations are derived as
$$FE_1=\left( \#f +n\# \nabla f\right)/n_1$$
and
$$FE_2=\left( \#f +n\# \nabla f+\frac{1}{2}n(n+1) \#\nabla^2f\right)/n_2,$$
where $\# f$, $\# \nabla f$, and $\# \nabla^2f$ denote the number of objective function, objective gradient, and objective Hessian evaluations. Performance profiles between MOSQP and MOSQP(F) using average function evaluations in best and worst run in \textit{\textbf{RAND}} are provided in Figures \ref{ctb} and \ref{ctw}. Figure 11 represents performance profiles between MOSQP and MOSQP(F) using average function evaluations in \textit{\textbf{LINE}}.
\begin{figure}[!htbp]
    \centering
    \begin{subfigure}[b]{\textwidth}
    \centering
     \includegraphics[height=2.5cm,width=.75\textwidth]{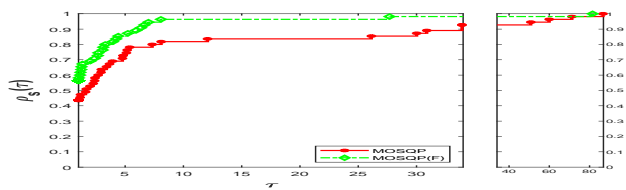}
     \label{ctl}
     \end{subfigure}
     \caption{Performance profile between MOSQP and MOSQP(F) using average function evaluations in \textit{\textbf{LINE}}}
\end{figure}
\\
{\bf Result analysis:} One may observe from the above figures that the method proposed in this article (MOSQP) gives better results than MOSQP(F) in purity, $\Gamma$, and $\Delta$ metrics using initial point selection strategy \textit{\textbf{RAND}} and purity, $\Gamma$ metrics using initial point selection strategy \textit{\textbf{LINE}}. Similarly MOSQP gives better results than MOS in $\Gamma$ metric using initial point selection strategy \textit{\textbf{RAND}} and \textit{\textbf{LINE}}. Other metrics have average performance ratios with MOS and MOSQP(F).
\section{Conclusion}
In this article we have developed a globally convergent modified SQP method for constrained multi-objective optimization problem. This method is free from any kind of a priori chosen parameters or ordering information of objective functions. Also feasibility of the sub-problem is guaranteed. To generate an approximate Pareto front, we have used the initial point selection strategies \textit{\textbf{LINE}} and \textit{\textbf{RAND}}. There is no single spreading technique for line search methods that can work in a satisfactory manner for all types of multi-objective programming problems. Spreading out an approximation to a Pareto front is a difficult task. A well distributed spreading technique is discussed in Step 3 of Algorithm 1.4 of \cite{flg3}. We keep the implementation of these techniques for future developments.

\end{document}